\documentclass[12pt,a4paper]{article}
\usepackage{amssymb,enumerate}
\usepackage{amsgen,amsmath,amstext,amsbsy,amsopn,amsfonts,amssymb}

\newtheorem{Thm}{Theorem}[section]
\newtheorem{Prop}[Thm]{Proposition}
\newtheorem{Lem}[Thm]{Lemma}
\newtheorem{Cor}[Thm]{Corollary}
\newtheorem{Def}[Thm]{Definition}
\newtheorem{Ex}[Thm]{Example}

\newtheorem{Rem}[Thm]{Remark}

\def\Ga{\Gamma}
\def\ga{\gamma}
\def\la{\lambda}
\def\H{\mathcal H}
\def\B{\mathcal B}
\def\K{\mathcal K}
\def\M{\mathcal M}

\def\NN{\mathbf N}
\def\CCC{\mathbf C}
\def\RRR{\mathbf R}
\def\QQ{\mathbf Q}

\newcommand{\bpr}{\noindent \textbf{Proof}: ~}

\newcommand{\epr}{~$\blacksquare$}

\newcommand{\eprsk}{~$\blacksquare$\medskip}

\title{Irreducible affine isometric actions on Hilbert spaces}
\author{Bachir Bekka, Thibault Pillon\footnote{Supported by grant 20-149261/1 of Swiss National Fund for scientific research} \;and Alain Valette}

\begin{document}
\maketitle
\begin{abstract}
Let $G$ be  locally compact group. We undertake a systematic study of irreducible  affine isometric actions of $G$  on Hilbert spaces.
It turns out that, while that are a few parallels of this study to the by now classical  theory of irreducible unitary representations,
these two theories differ in several aspects (for instance, the direct sum of two irreducible affine actions can still be irreducible).
One of the main tools we use is an affine  version of Schur's  lemma characterizing  the irreducibility of an affine isometric action of $G$.
This enables us to describe for instance  the irreducible  affine isometric actions of  nilpotent groups.
As another application,  a short proof is provided  for the following result of Neretin:  the restriction  to a cocompact lattice  of an irreducible affine action of $G$   remains irreducible.  We give a necessary and sufficient condition for a fixed unitary representation $\pi$  to be
the linear part of an irreducible affine action. In particular, when $\pi$ is a multiple of the regular representation of a discrete group $\Gamma$,  we show how this question is related to the  $L^2$-Betti number $\beta^1_{(2)}(\Gamma).$ 
After giving  a necessary and sufficient condition for a direct sum of irreducible affine actions to be irreducible,  we show the following super-rigidity result: if $G$ is product of two or more locally compact groups and $\Gamma$  an irreducible co-compact lattice in $G,$ then  any irreducible affine action $\alpha$ of $\Ga$ extends to an affine action of $G$, provided the linear part of $\alpha$ does not weakly contain the trivial representation.

\end{abstract}
\section{Introduction}

The theory of unitary representations of locally compact groups is by now a central and classical part of representation theory. Very quickly, the theory centers on the study of unitary irreducible representations which, for suitable classes of groups  (e.g. compact Lie groups, nilpotent Lie groups, semi-simple Lie groups, to name just a few), has reached a very satisfactory state. 

The theory of affine isometric actions on Hilbert spaces is, comparatively, a much more recent subject, that developed through connections with property (T), the Haagerup property, or operator algebras (see e.g. \cite{BHV}). To the best of our knowledge, {\it irreducible} affine isometric actions were first considered by Neretin \cite{Ner}, who also provides many examples. So let $\alpha$ be an affine isometric action of the group $G$ on the Hilbert space ${\cal H}$, i.e. a group homomorphism $\alpha:G\rightarrow {\rm Isom}({\cal H})$ from $G$ to the group of affine isometries  of $\cal H$.

\begin{Def}  \rm The action $\alpha$ is {\it irreducible} if ${\cal H}$ has no non-empty, closed and  proper $\alpha(G)$-invariant affine subspace.
\end{Def}

The following two classes of examples should be kept in mind.

\begin{Ex} \rm Let $b:G\rightarrow{\cal H}$ be a homomorphism to the additive group of ${\cal H}$. It gives rise to an affine action of $G$ by translations on ${\cal H}$, which is irreducible if and only if the linear span of $b(G)$ is dense in ${\cal H}$.
\end{Ex}

\begin{Ex} \rm\label{basicex2} Let $\pi$ be a unitary irreducible representation of $G$ on ${\cal H}$, such that $H^1(G,\pi)\neq 0$. Choose a 1-cocycle $b$ which is not a 1-coboundary. Then the affine action $\alpha$ of $G$ on ${\cal H}$, defined 
$$\alpha(g)v:=\pi(g)v+b(g) \qquad\text{for} \qquad g\in G,\, v\in{\cal H},$$
 is irreducible. Indeed, assume by contradiction that ${\cal K}$ is a non-empty, closed, proper, $\alpha(G)$-invariant affine subspace. Then its linear part ${\cal K}_0,$ 
 is  a proper and closed $\pi(G)$-invariant linear subspace;  by irreducibility of $\pi,$ it follows that ${\cal K}_0=0$. So $\alpha$ has a fixed point, contradicting the fact that $b$ is not a coboundary.
\end{Ex}

In this paper, we undertake a systematic study of irreducible  affine isometric actions of the locally compact group $G$  on Hilbert spaces.
The theory of irreducible affine isometric actions has some parallels with the theory of irreducible unitary representations, but to a limited extent. To illustrate this, we contrast the classical case and the affine case in two columns (where the left column is about a unitary representation $\pi$, the right column is about an affine action $\alpha$ with linear part $\pi$ and translation part $b$).

\begin{enumerate}
\item Characterization

\begin{tabular}{c|c}
\begin{minipage}{75mm}
$\pi$ is irreducible if and only $\pi(G)\xi$ is total for every non-zero vector $\xi$ if and only if every positive-definite function $g\mapsto\langle\pi(g)\xi|\xi\rangle$ lies on an extremal ray in the cone of positive-definite functions on $G$.
\end{minipage} & 
\begin{minipage}{75mm} $\alpha$ is irreducible if and only if, for every vector $v$, the cocycle $g\mapsto b(g)+\pi(g)v-v$ has total image; if and only if $b(G)$ is total and, for every decomposition $\|b(g)\|^2=\psi_0(g)+\psi_1(g)$, with $\psi_0,\psi_1$ functions conditionally of negative type with $\psi_0\neq 0$, the function $\psi_0$ is unbounded (see Proposition \ref{affine}).

\end{minipage}
  \end{tabular}

\item Existence ($G$ locally compact)

  \begin{tabular}{c|c}
\begin{minipage}{75mm}
Irreducible  unitary representations of $G$ separate points (Gelfand-Raikov).
\end{minipage} & 
\begin{minipage}{75mm} For $G$ compactly generated, $G$ admits an irreducible affine action if and only if $G$ doesn't have property (T), as follows from Theorem~0.2 in \cite{ShaInv}. Even then, irreducible affine actions do not separate points in general (see Corollary \ref{nilp} below).
\end{minipage}
  \end{tabular}


\item Commutants

 \begin{tabular}{c|c}
\begin{minipage}{75mm}
$\pi(G)'$ is the commutant of $\pi(G)$ in $B(\H)$ (it is a von Neumann algebra)
\end{minipage} & 
\begin{minipage}{75mm} $\alpha(G)'$ is the commutant of $\alpha(G)$ in the monoid of continuous affine maps on $\H$. The affine map $Av=:Tv+t$ is in $\alpha(G)'$ iff $T\in\pi(G)'$ and $(T-1)b(g)=\pi(g)t-t$ for all $g\in G$ (see Lemma \ref{easylem}).
\end{minipage}
  \end{tabular}

\item Schur's lemma
 
  \begin{tabular}{c|c}
\begin{minipage}{75mm}
$\pi$ is irreducible iff $\pi(G)'=\CCC.1$
\end{minipage} & 
\begin{minipage}{75mm} $\alpha$ is irreducible if and only if $\alpha(G)'$ consists of translations (in this case, exactly the set of translations along $\H^{\pi(G)}$; see Proposition \ref{Schur}).
\end{minipage}
  \end{tabular}


\item Abelian groups
 
  \begin{tabular}{c|c}
\begin{minipage}{75mm}
Every  irreducible unitary representation is one-dimensional
\end{minipage} & 
\begin{minipage}{75mm} Every irreducible action is given by some homomorphism $b:G\rightarrow\H$ with $b(G)$ having dense linear span (see Proposition \ref{FC-irr}).
\end{minipage}
  \end{tabular}   
   
\item Nilpotent groups
   
    \begin{tabular}{c|c}
\begin{minipage}{75mm}
Usually, the irreducible unitary representations of $G$ are infinite dimensional (think of Kirillov's orbit method).
\end{minipage} & 
\begin{minipage}{75mm} Same as for abelian groups, see Corollary \ref{nilp}.
\end{minipage}
  \end{tabular}
  
  \end{enumerate}

Apart from allowing us to determine the irreducible affine actions of abelian or nilpotent groups, our affine Schur lemma has several other applications:
\begin{itemize}
\item We give in Theorem \ref{neretin} a short proof of Neretin's result \cite{Ner} that, upon restricting to a co-compact lattice in a locally compact group, an irreducible affine action remains irreducible\footnote{It is well-known that, in general, restricting a unitary irreducible representation to a co-compact lattice, does not yield an irreducible representation.}.
\item We are able to study the question: ``when is a given unitary representation $\pi$ the linear part of an irreducible affine action?'' In particular, taking for $\pi$ a multiple of the regular representation of a non-amenable, ICC discrete group $\Gamma$, we get a new definition of the first $L^2$-Betti number $\beta^1_{(2)}(\Gamma)$; namely $\beta^1_{(2)}(\Gamma)$ is the supremum of all non-negative $t$'s such that the unique module over the von Neumann algebra $L(\Gamma)$ of $\Gamma$ with $L(\Gamma)$-dimension $t$ is the linear part of some irreducible affine action (see Corollary \ref{ICCnewbetti}).
\item The definition of $L^2$-Betti numbers $\beta^n_{(2)}$ has been extended from discrete to locally compact unimodular groups, in two papers by Petersen \cite{Pet} and Kyed-Petersen-Vaes \cite{KPV}. We prove in Theorem \ref{betti+discreteser} that, if $G$ is a locally compact group containing a co-compact lattice, then 
\begin{equation}\label{inegal}
\beta^1_{(2)}(G)\geq \sum_{\sigma\in\hat{G}_d} d_\sigma\dim_\CCC H^1(G,\sigma)
\end{equation}
where $\hat{G}_d$ is the discrete series of $G$ (i.e. the set of square-integrable unitary irreducible representations of $G$, up to unitary equivalence), and $d_\sigma>0$ is the formal dimension of  $\sigma$. The proof depends crucially on irreducible affine actions, even if the inequality involves no such actions .
\end{itemize}
Here is a short summary of the paper.  We give in Section 2 a number of characterizations of irreducible affine actions. Commutants are introduced in Section 3, where the affine Schur lemma is also proved. Section 4 contains several applications of the affine Schur lemma: to restricting affine actions to lattices, to the behavior of an irreducible affine action on the center of a group, to abelian and nilpotent groups, and to the regular representation of a discrete group. Observing that (unlike what happens for unitary representations!), the direct sum of two irreducible affine actions can still be irreducible, we give in Section 5 a necessary and sufficient condition for this to happen. In Section 6, we combine this with a super-rigidity result of Shalom \cite{ShaInv} and show that, if $\Gamma$ is an irreducible co-compact lattice in a product of two or more locally compact groups, any irreducible affine action of $\Gamma$ extends to an affine action of the ambient group, provided the linear part of $\alpha$ does not weakly contain the trivial representation. Section 7 is devoted to the proof of inequality (\ref{inegal}) mentioned above. Finally in Section 8 we compare our notion of irreducibility for affine actions with other possible definitions, already introduced in \cite{CTV-SO}.

\medskip
\noindent
{\bf Ackowledgements:} We thank P.-E. Caprace, T. Gelander, N. Monod, A. Thom, S. Vaes for useful conversations at various stages of the project.

\section{Characterizations of irreducible affine actions}

\subsection{Notations}\label{not}

Let $G$ be a topological group with identity element $e$; a continuous function $\psi:G\rightarrow\mathbf{R}$ is {\it conditionally of negative type} (CNT) if $\psi(e)=0$ and, for every $n\geq 1,\,g_1,\dots,g_n\in G$, and $\lambda_1,\dots,\lambda_n\in\mathbf{R}$ with $\sum_{i=1}^n\lambda_i=0$, we have 
$$\sum_{i,j=1}^n \lambda_i\lambda_j\psi(g_i^{-1}g_j)\leq 0.$$
Equivalently, by the GNS construction (see \cite{BHV}, Theorem C.2.3), there exists a Hilbert space ${\cal H}_\psi$ and a (continuous) affine isometric action $\alpha_\psi$ of $G$ on ${\cal H}_\psi$ such that $\psi(g)=\|\alpha_\psi(g)(0)\|^2$ for every $g\in G$.

Let $C$ be the cone of CNT functions on $G$. It is known (see \cite{VK}, or Th\'eor\`eme 1 in \cite{LSV} \footnote{Note that the assumption $b\neq 0$ is missing in the statement of this result in \cite{LSV}; also, it should have been said in the proof that the linear subspace spanned by $b(G)$ is $\pi(G)$-invariant (as follows easily from the 1-cocycle relation), hence by irreducibility it is dense in ${\cal H}$.}) that a non-zero $\psi\in C$ lies on an extremal ray, if and only if the linear part $\pi_\psi$ of $\alpha_\psi$ is an irreducible orthogonal representation of $G$. Define two sub-cones $C_b$ and $C_u$: the cone $C_b$ is the set of bounded functions in $C$, and the cone $C_u$ is the set of unbounded functions in $C$, together with $0$. Clearly $C=C_b\cup C_u$, and $C_b\cap C_u=\{0\}$, and $C_b$ is a face in $C$. For $G$ locally compact $\sigma$-compact group, $C_u=\{0\}$ if and only if $G$ has Kazhdan's property (T): this is a re-phrasing of the Delorme-Guichardet theorem (see \cite{BHV}, Theorem 2.12.4).

 Let $(\pi,{\cal H})$ be a unitary representation of $G$ on a Hilbert space ${\cal H}$; we denote by $Z^1(G,\pi)$ (resp. $B^1(G,\pi)$) the space of 1-cocycles (resp. 1-coboundaries) associated with $\pi$. The 1-cohomology $H^1(G,\pi)$ is the quotient $Z^1(G,\pi)/B^1(G,\pi)$.
 
Let $b\in Z^1(G,\pi)$ be a 1-cocycle. Let $\psi(.)=\|b(.)\|^2$ be the associated function conditionally of negative type, and $\alpha_{\pi,b}$ the associated affine isometric action of $G$ on ${\cal H}$, defined by $\alpha_{\pi,b}(g)v=\pi(g)v+b(g)$ (for $g\in G,v\in {\cal H})$. When $\pi$ and $b$ are clear, we will write $\alpha$ for $\alpha_{\pi,b}$.

\subsection{Characterizations of irreducibility}\label{CharIrr}

For $v\in{\cal H}$, we shall denote by $\partial_v$ the 1-coboundary $\partial_v(.):=\pi(.)v-v$; this is the 1-cocycle associated with the affine isometric action $t_v^{-1}\circ\pi\circ t_v$, where $t_v$ is the translation of vector $v$ in ${\cal H}$, so this affine action has a fixed point and it is reducible. Thus {\it we will assume from now on that $b$ is not a 1-coboundary.} Throughout, all affine subspaces will be assumed to be non-empty.

Let $\pi_0$ be a sub-representation of $\pi$, on a closed subspace $V_0\subset{\cal H}$. Let us denote by $b_0(g)$ the orthogonal projection of $b(g)$ on $V_0$. It is immediate to check that $g\mapsto b_0(g)$ is a cocycle with respect to $\pi_0$, so that $\alpha_0(g)v=\pi_0(g)v+b_0(g)$ defines an affine isometric action of $G$ on $V_0$: we call it the {\it projected action} on $V_0$.

Recall that a subset of ${\cal H}$ is {\it total} if it generates a dense linear subspace of ${\cal H}$.

\begin{Prop}\label{affine} 
Keep notations as in subsection \ref{not}. The following properties are equivalent:
\begin{enumerate}
\item[(A1)] The affine isometric action $\alpha$ is irreducible.
\item[(A2)] For every $v\in{\cal H}$, the 1-cocycle $b+\partial_v$ has total image in ${\cal H}$.
\item[(A3)] For every direct sum decomposition $\pi=\pi_0\oplus\pi_1$ with $\pi_0\neq 0$, in the corresponding decomposition $b=b_0\oplus b_1$, the 1-cocycle $b_0$ is unbounded.
\item[(A4)] $b(G)$ is total and, for every decomposition $\psi=\psi_0+\psi_1$, with $\psi_0,\psi_1$ functions conditionally of negative type with $\psi_0\neq 0$, the function $\psi_0$ is unbounded.
\item[(A5)] $b(G)$ is total and $\psi$ belongs to a common face of $C$ and $C_u$.
\item[(A6)] For every non-zero sub-representation $\pi_0$ of $\pi$, the projected action $\alpha_0$ is irreducible.
\end{enumerate}
\end{Prop}

\bpr We follow the schemes $(A1)\Rightarrow(A6)\Rightarrow(A3)\Rightarrow(A2)\Rightarrow(A1)$ and $(A1)\Rightarrow(A4)\Leftrightarrow(A5)\Rightarrow(A3)$

$(A1)\Rightarrow(A6)$ Assume that there is a closed, $\pi(G)$-invariant subspace $V_0\subset{\cal H}$ such that the projected action $\alpha_0$ is reducible. So there exists a proper closed, $\alpha_0(G)$-invariant affine subspace $W\subset V_0$. Let $V_0^\perp$ denote the orthogonal complement of $V_0$. Then $W\oplus V_0^\perp$ is a proper closed, $\alpha(G)$-invariant affine subspace of ${\cal H}$, so that $\alpha$ is reducible.

$(A6)\Rightarrow(A3)$ is clear, as boundedness of $b_0$ implies reducibility of $\alpha_0$.

$(A3)\Rightarrow(A2)$: Assume that, for some $v\in{\cal H}$, the set $(b+\partial_v)(G)$ is not total. Let $W_1$ be the closed linear subspace it generates. It follows from the 1-cocycle relation for $b+\partial_v$ that $W_1$ is $\pi(G)$-invariant. Let $W_0$ be the orthogonal complement of $W_1$, and let $\pi=\pi_0\oplus\pi_1, \; b=b_0\oplus b_1$ and $v=v_0\oplus v_1$ be the corresponding decompositions of $\pi,\;b$, and $v$. As $v+W_1$ is $\alpha(G)$-invariant, it follows that the affine action $\alpha_0$ obtained by projecting to $W_0$ has $v_0$ as a fixed point, i.e. $b_0$ is bounded.

$(A2)\Rightarrow(A1)$: Assume by contraposition that $\alpha$ has a non-empty, closed invariant affine subspace $W\neq{\cal H}$; let $W_0=W-W$ be the corresponding linear subspace, so that $W_0\neq{\cal H}$. Then for $v\in W$ we have $\alpha(g)v-v\in W_0$ for every $g\in G$, i.e. $b(g) + \pi(g)v-v\in W_0$, showing that $(b+\partial_v)(G)$ is not total.

 $(A1)\Rightarrow(A4)$ We proceed by contraposition. If $b(G)$ is not total, then $\alpha$ is reducible. Suppose now that $b(G)$ is total and there exists a decomposition $\psi=\psi_0+\psi_1$ where $\psi_0$ is non-zero and bounded. Let $\alpha_0$ be the affine isometric action associated with $\psi_0$ by the GNS construction (see Proposition 2.10.2 in \cite{BHV}). It has a fixed point $w$, as $\psi_0$ is bounded. Now, by the proof of Theorem 1 in \cite{LSV} (see in particular pp. 245-246) the map $\sum_i a_ib(g_i)\mapsto \sum_i a_ib_0(g_i)$, from the span of $b(G)$ to the span of $b_0(G)$, extends linearly and continuously to a bounded linear map $T_0:{\cal H}\rightarrow{\cal H}_0$, which is onto and intertwines $\alpha$ and $\alpha_0$. Hence $T_0^{-1}(w)$ is a proper, closed, affine subspace of ${\cal H}$ which is $\alpha(G)$-invariant, so $\alpha$ is reducible.
 
$(A4)\Rightarrow (A5)$: Set $F=\{\psi_0\in C: \exists \psi_1\in C\,\mbox{such that}\,\psi_0+\psi_1\in\mathbf{R}^+\psi\}$. This is clearly the smallest face of $C$ containing $\psi$. The assumption implies that $F\subset C_u$, so $F$ is a common face of $C$ and $C_u$.
 
$(A5)\Rightarrow (A4)$ is obvious.
 
$(A4)\Rightarrow (A3)$: Set $\psi_i(.)=\|b_i(.)\|^2\;(i=0,1)$ and notice that the assumption that $b(G)$ is total implies that $b_0\neq 0$.
\eprsk

\begin{Ex} \rm If $\alpha$ is irreducible then by $(A1)\Rightarrow(A2)$ the set $b(G)$ is total in ${\cal H}$. The converse is {\it false}: the reason is that condition $(A2)$ is translation-invariant, while $b(G)$  being total is not. Concretely, let $G=\mathbf{Z}$ act isometrically on $\mathbf{R}^2$ by 
$$
\alpha_n(x,y)=(x+n,(-1)^ny+1-(-1)^n)\quad \text{for all} \quad n\in\mathbf{Z},\;(x,y)\in\mathbf{R}^2.
$$ Geometrically, this is the action by powers of the glide symmetry with axis the horizontal line $y=1$, and translation by $+1$ to the right. Then \emph{all} orbits are total, in particular $\alpha(G)(0)=b(G)$, but $\alpha$ is reducible as the axis is invariant.
\end{Ex}

\section{Use of commutants}

\subsection{The commutant of an affine action}

Let $\alpha$ be an affine isometric action of a group $G$, with linear part $\pi$. We recall that the {\it commutant} of $\pi$ is the von Neumann algebra 
$$\pi(G)'=\{T\in{\cal B(H)}: T\pi(g)=\pi(g)T \quad \text{for all} \quad  g \in G\}.$$
 If $b$ is a cocycle for $\pi$ and $T\in\pi(G)'$, we observe that $Tb$ is still a cocycle for $\pi$, so that $\pi(G)'$ acts on the space $Z^1(G,\pi)$ of 1-cocycles, and this action descends to the first cohomology space $H^1(G,\pi)$.

\begin{Def} \rm  
The {\it commutant} of $\alpha$ is the set of (continuous) affine transformations $A$ on ${\cal H}$ such that $A\circ\alpha(g)=\alpha(g)\circ A$ for every $g\in G$.
\end{Def}

Write an affine transformation $A$ on ${\cal H}$ as $Av=Tv+t$ for $v\in{\cal H}$, where $T\in{\cal B(H)}$ is the linear part. It is easy to see that $A$ is in the commutant of $\alpha$ if and only if $T\in\pi(G)'$ and $(T-1)b(g)=\partial_t(g)$ for every $g\in G$. From this the following lemma is immediate:

\begin{Lem}\label{easylem} For $T\in\pi(G)'$, the following properties are equivalent:
\begin{enumerate}
\item[i)] There exists $t\in{\cal H}$ such that the affine transformation $Av=:Tv+t$ is in the commutant of $\alpha$.
\item[ii)] There exists $t\in{\cal H}$ such that $(T-1)b(g)=\partial_t(g)$ for every $g\in G$.
\item[iii)] $(T-1)[b]=0$, where $[b]$ denotes the class of $b$ in $H^1(G,\pi)$. \epr
\end{enumerate}
\end{Lem}

\begin{Rem}
\rm
 We observe that, if $Av=Tv+t$ is in the commutant of an affine action $\alpha$ without fixed point, then 1 is a spectral value of $T$, as the operator $T-1$ maps the unbounded set $b(G)$ to the bounded set $\partial_t(G)$.
\end{Rem}

\subsection{A Schur-type lemma}

We denote by ${\cal H}^{\pi(G)}$ the space of $\pi(G)$-fixed vectors in ${\cal H}$.

\begin{Prop}\label{Schur} Let $\alpha$ be an affine isometric action on ${\cal H}$. The following properties are equivalent.
\begin{enumerate}
\item[i)] $\alpha$ is irreducible.
\item[ii)] The commutant of $\alpha$ is the set of translations along ${\cal H}^{\pi(G)}$.
\item[iii)] The commutant of $\alpha$ consists of translations.
\end{enumerate}
\end{Prop}

\bpr $(i)\Rightarrow (ii)$ Let $Av=Tv+t$ be an affine transformation of ${\cal H}$, in the commutant of $\alpha$. Then $T\in\pi(G)'$ and 
\begin{equation}\label{commut}
(T-1)b(g)=\pi(g)t-t \quad \text{for every}\quad g\in G.
\end{equation}
So it is enough to show that $T=1$. For this, consider the positive operator 
$$S=T^*T-T-T^*+2 = (T-1)^*(T-1)+1;$$
 if we show $S=1$, then $T=1$. As $S$ is self-adjoint, it is enough to show that the spectrum of $S$ is $\{1\}$. Assume by contradiction that there some other spectral value $s$. Let $[a,b]$ be a closed interval of $\mathbf{R}$ containing $s$ in its interior, and not containing 1. Let $E$ be the spectral projector of $S$ associated with $[a,b]$, so that $E\neq 0$ and $E\in\pi(G)'$. Denote by $\rho$ the sub-representation of $\pi$ on ${\rm Im}(E)$. Apply $(T-1)^*$ to Equation \ref{commut}:
$$(S-1)b(g)=(\pi(g)-1)(T^*-1)t.$$
Then apply $E$ and restrict to ${\rm Im}(E)$:
$$(S-1)Eb(g) = (\rho(g)-1)E(T^*-1)t.$$
But $S-1$ is invertible as a bounded operator on ${\rm Im}(E)$ (since $1\notin [a,b]$); denoting by $R$ its inverse, we obtain
$$Eb(g)=(\rho(g)-1)RE(T^*-1)t.$$
The projection $Eb$ of $b$ on ${\rm Im}(E)$ is therefore bounded, contradicting condition (A3) in Proposition \ref{affine}.

$(ii)\Rightarrow(iii)$ is trivial.

$(iii)\Rightarrow(i)$ Assume that $\alpha$ is reducible, and let $W$ be a non-trivial closed, invariant, affine subspace of ${\cal H}$. Let $E:{\cal H} \rightarrow W$ be the projection onto $W$; so   $Ev$ is the point of $W$ closest to $v$, for every $v\in \cal H$.  
Since every $\alpha (g)$ is an isometry, it follows that   the affine transformation  $E$ is in the commutant of $\alpha$. \epr

\medskip
We already observed that the first cohomology $H^1(G,\pi)$ is a module over the von Neumann algebra $M:=\pi(G)'$; recall that a vector $\xi$ in a module over $M$, is {\it separating} if $S\xi=0$ implies $S=0$ for every $S\in M$.

\begin{Cor}\label{separating} Let $\pi$ be a unitary representation of $G$. There exists an irreducible affine action $\alpha$ with linear part $\pi$ if and only if $H^1(G,\pi)$ admits a separating vector for $\pi(G)'$.
\end{Cor}

\bpr According to Proposition \ref{Schur}, the existence of $\alpha$ is equivalent to the existence of a 1-cocycle $b$ such that, for every $T\in\pi(G)'$ and $t\in{\cal H}$ such that $(T-1)b(g)=\partial_t(g)$ for every $g\in G$, we have $T=1$; in turn, by Lemma \ref{easylem}, this is equivalent to the existence of a class $[b]\in H^1(G,\pi)$ such that $(T-1)[b]=0$ for $T\in\pi(G)'$, implies $T=1$; this exactly means  that $[b]$ is a separating vector for $\pi(G)'$.\epr

\section{Applications}

\subsection{Restriction to lattices}

We give a short proof of a result of Neretin (Theorem 3.6 in \cite{Ner}\footnote{We seize this opportunity to correct an error in \cite{Ner}: the proof of Theorem 3.6 rests on Proposition 2.5 of the same paper, which claims that, if an affine isometric action $\alpha$ has a closed, affine invariant subspace $L$ such that $\alpha|_L$ is irreducible, then every closed, affine invariant subspace of $\alpha$ contains $L$: this is false, as shown by an action of $\mathbf{Z}$ by translations on the plane. It can be checked however that Neretin's proof holds for irreducible affine actions whose linear part has no non-zero fixed vector.}) asserting that the restriction of an irreducible affine action to a co-compact lattice, remains irreducible. Since we do not use induction of affine actions, we are able to remove the assumption of discreteness of the subgroup in \cite{Ner}. In order to treat non-co-compact lattices, we introduce a definition: for $H$ a lattice in $G$ and $b\in Z^1(G,\pi)$, we say that the cocycle $b$ is { \it integrable on $G/H$} if there exists a measurable fundamental domain $\Omega$ for the right action of $H$ on $G$, such that $\int_\Omega\|b(g)\|\,dg<+\infty$, where $dg$ denotes Haar measure on $G$.

\begin{Thm}\label{neretin} Let $H$ be a closed subgroup of the locally compact group $G$, such that $G/H$ carries a $G$-invariant probability measure $\mu$. Let $\alpha(g)v=\pi(g)v+b(g)$ be an affine isometric action of $G$. Assume  either that $H$ is co-compact or that $H$ is discrete and the cocycle $b$ is integrable on $G/H$. If $\alpha$ is irreducible, then the restriction $\alpha|_H$ is irreducible.
\end{Thm} 

\bpr Let ${\cal K}$ be a closed affine subspace, invariant under $\alpha|_H$, and let $E$ be the projection onto ${\cal K}$. We want to show that $E$ is the identity of ${\cal H}$, or equivalently that its linear part $E_0$ is the identity. Write $Ev=E_0v + t$ for $v\in{\cal H}$

Let $\textnormal{Aff}({\cal H})$ be the set of continuous affine maps from ${\cal H}$ to ${\cal H}$. Consider the map 
$$G\rightarrow\textnormal{Aff}({\cal H}):g\mapsto \alpha(g)E\alpha(g)^{-1};$$
 this map factors through $G/H$, and we wish to integrate it on $G/H$. For this, we compute (using $b(g^{-1})=-\pi(g)^{-1}b(g)$):
$$\alpha(g)E\alpha(g)^{-1}v=\pi(g)E_0\pi(g)^{-1}v + \pi(g)t+[1-\pi(g)E_0\pi(g)^{-1}]b(g).$$
The first two terms are bounded, and the third one is integrable on $G/H$ under either of our assumptions. So we may define 
\begin{equation}\label{average}
A=\int_{G/H} \alpha(x)E\alpha(x)^{-1}\,d\mu(x)
\end{equation}
as an element of $\textnormal{Aff}({\cal H})$. By $G$-invariance of $\mu$, we see that $A$ belongs to the commutant of $\alpha$. By Proposition \ref{Schur}, the affine transformation $A$ is a translation. Taking linear parts in Equation (\ref{average}), we get $1=\int_{G/H}\pi(x)E_0\pi(x)^{-1}\,d\mu(x)$, expressing the identity $1$ on $\cal H$ as an average of operators of norm $\leq 1$. Since $1$ is an extreme point in the unit ball of ${\cal B(H)}$ (see e.g. Proposition 1.4.7 in \cite{Ped}), we deduce $E_0=1$.\epr

\medskip

\begin{Rem} 
\label{Rem-IntegrableCocycle}
\rm
Let us take a closer look at the condition of integrability of the cocycle in the case of a non-uniform lattice $\Gamma$ in $G$. Assume that the ambient group $G$ is compactly generated, and denote by $|g|_S$ the word length of $g\in G$ with respect to some compact generating set $S\subset G$. If $b\in Z^1(G,\pi)$, it is an easy consequence of the triangle inequality that there exists $C>0$ such that $\|b(g)\|\leq C|g|_S$;
 so,  for  a lattice $\Gamma$ in $G$,  a sufficient condition for every cocycle to be integrable on $G/\Gamma$ is the existence of a measurable fundamental domain $\Omega$  for the right action of  $\Gamma$ on $G$ such that:
\begin{equation}\label{funddomain}
\int_\Omega |g|_S\,dg<\infty
\end{equation}

This is of course clear for uniform lattices. Margulis proves it for S-arithmetic groups in \cite[Prop. VIII.1.2]{Mar}. Using the Garland-Raghunathan description of cusps \cite{GR}, it can be checked that this condition is also satisfied by all lattices in rank 1 simple Lie groups. It also holds for twin buildings lattices, see \cite[Lemma 4.2]{CapRem}.

It was however pointed out to us by Tchachik Gelander that the condition in (\ref{funddomain}) does not hold in general, as counterexamples can be found into the automorphism group ${\rm Aut}(T_k)$ of the $k$-regular tree, with $k\geq 3$. Indeed, consider the graph of groups based on the infinite ray with vertices $x_0,x_1,x_2,\dots$. Denote by $\Gamma_n$ the vertex group at $x_n$, and $H_n$ the edge group at the edge $[x_{n-1},x_n]$ (for $n\geq 1$). Assume that indices satisfy $$[\Gamma_n:H_n]+[\Gamma_n:H_{n+1}]=k,$$ so that the fundamental group $\Gamma$ of the graph of groups (in the sense of Bass-Serre \cite{Ser}) acts on $T_k$. Assume now that the $\Gamma_n$'s are finite groups, whose orders satisfy 
$$\sum_{n=1}^\infty\frac{1}{|\Gamma_n|}<+\infty \qquad \text{but} \qquad \sum_{n=1}^\infty\frac{n}{|\Gamma_n|}=+\infty.$$
 The former condition ensures that $\Gamma$ sits in ${\rm Aut}(T_k)$ as a non-uniform lattice (see \cite{Ser}, Section 1.5 in Chapter II), while the latter condition implies the non-existence of $\Omega\subset {\rm Aut}(T_k)$ such that (\ref{funddomain}) holds. The construction of the $\Gamma_n$'s requires some care, due to the constraints on the indices of $H_n$ and $H_{n+1}$. For example, one can define a sequence $(a_i)_{i\geq 0}$ of positive integers in a recursive way, by requiring 
 $$a_0=1\qquad \text{and} \qquad  a_{i}-a_{i-1}=\lfloor\frac{(k-1)^i}{i^2}\rfloor,$$
  and then choose $H_n$ with $|H_n|=(k-1)^i$ for $a_i\leq n<a_{i+1}$.
\end{Rem} 

\begin{Rem} 
\label{Rem-Induction}
\rm
Let $\Gamma$ be a co-compact lattice in the locally compact group $G$. Given   an action $\alpha$ of $\Ga$ by affine isometries
on a Hilbert space $\cal H,$   it is possible to  define
an \emph{induced} affine action ${\rm Ind}_\Ga^G \alpha$  of $G$, as  discussed in  \cite[Section II]{ShaInv}.
Let us briefly review the construction.
Let  $\pi$ be the  linear part  of $\alpha$  and  $b\in Z^1(\Gamma, \pi)$ the corresponding $1$-cocycle.
Let $\Omega$ be a compact  fundamental domain for the right action of $\Gamma$ on $G$
and  $c: G\times \Omega\to \Gamma$   the associated cocycle  defined by $c(g,x) =\gamma$ if and only if 
$gx\gamma \in \Omega.$ 
The induced unitary representation ${\rm Ind}_\Ga^G \pi$ of $G$ can be realized on 
$L^2(\Omega, \H)$ by means of the formula
$$
({\rm Ind}_\Ga^G \pi)(g) f(x) = \pi(c(g^{-1}, x)) f(g^{-1}x) \qquad f\in L^2(\Omega, \H), \, g\in G, \, x\in \Omega.
$$
The map  $\widetilde{b}: G\to L^2(\Omega, \H),$ defined by 
$$
\widetilde{b}(g) (x)= b(c(g^{-1}, x)) \qquad \, g\in G, \, x\in \Omega,
$$
belongs to $Z^1(G, {\rm Ind}_\Ga^G \pi)$; observe that, since $\Omega$ is compact,  $\widetilde{b}$ takes
indeed its values in $L^2(\Omega, \H).$ The induced affine action ${\rm Ind}_\Ga^G\alpha$ of $G$ is  the action with linear  
part ${\rm Ind}_\Ga^G\pi$ and translation part given by $\widetilde{b}.$

One may ask whether ${\rm Ind}_\Ga^G\alpha$ is irreducible when $\alpha$ is irreducible.
This is not the case, even when $\Gamma$ has finite index in $G,$ as the following simple example shows.
Let $G= C_2\times \mathbf Z$ be the direct product of the cyclic group   of order two and the group of integers
and let $\Gamma= \mathbf Z$. Let $\alpha$ be the affine isometric action of $\Gamma$ on $\mathbf R$ 
defined  by 
$$
\alpha(n)y= y+n, \qquad n\in \mathbf Z,\,  y\in \mathbf R.
$$
So, the linear part of $\alpha$ is the identity and the injection $\mathbf Z\to \mathbf R$
is the corresponding cocycle.
The induced affine action  ${\rm Ind}_\Ga^G\alpha$ of $G$ is  easily seen to be defined 
on ${\mathbf R}^2$ by 
$$
({\rm Ind}_\Ga^G\alpha)(a,n)(x) = (x, y+n) \qquad n\in \mathbf Z, \, (x,y)\in {\mathbf R}^2.
$$
Clearly,  ${\rm Ind}_\Ga^G\alpha$ is not irreducible.

\end{Rem}

\subsection{Center and FC-center}

We denote by $Z(G)$ the center of the topological group $G$.

\begin{Prop}\label{center} In an irreducible affine action $\alpha$ of $G$ on ${\cal H}$, the center $Z(G)$ acts by translations in the direction of ${\cal H}^{\pi(G)}$. 
\end{Prop}

\bpr This follows immediately from Proposition \ref{Schur}. \epr

\begin{Cor}\label{G/Z} Assume that $\textnormal{Hom}(G,\mathbf{R})=0$. Then every irreducible affine action $\alpha$ of $G$ factors through $G/Z(G)$.
\end{Cor}

\bpr Let $b$ be the cocycle defining $\alpha$, and let $b_0$ be its projection on ${\cal H}^{\pi(G)}$, so that $b_0$ is a continuous homomorphism from $G$ to the additive group of ${\cal H}^{\pi(G)}$, hence $b_0\simeq 0$ by our assumption. This forces ${\cal H}^{\pi(G)}=0$ (otherwise we would contradict condition (A3) in Proposition \ref{affine}). By Proposition \ref{center}, the center $Z(G)$ acts by the identity. \epr

As a consequence, we get a very short proof of a result of J.-P. Serre (see Theorem 1.7.11 in \cite{BHV}).

\begin{Cor}\label{Serre} Let $G$ be a compactly generated, locally compact group. Assume that the separated abelianization $G/\overline{[G,G]}$ is compact. Let $Z$ be a closed central subgroup of $G$. If $G/Z$ has property (T), then so does $G$.
\end{Cor}

\bpr Our assumption implies that ${\rm Hom}(G,\mathbf{R})=0$. Assume by contraposition that $G$ does not have property (T). Since $G$ is compactly generated, the group $G$ admits an irreducible affine action $\alpha$,  by Shalom's theorem (\cite[Theorem 0.2]{ShaInv}). By Corollary \ref{G/Z}, this action $\alpha$ is actually an irreducible affine action of $G/Z$, which therefore does not have property (T). \epr

\medskip
The {\it FC-center} of $G$, denoted  ${\rm FC}(G),$  is  the set of  elements in $G$ 
with finite conjugacy class. Observe that the conjugacy class
of an element $\ga$ is finite if and only its
 centralizer $C_\ga$  in $G$ has finite index in $G.$
The FC-center  is a  subgroup of $G$
which is of course characteristic.

Observe that the FC-center of any group $\Ga$ is amenable. 
Indeed, every finitely generated subgroup of ${\rm FC}(\Ga)$
has a center of finite index and is hence amenable; it follows that ${\rm FC}(\Ga)$
is a union of amenable groups and is therefore amenable.

\begin{Prop}\label{FC-center} Let $\alpha$ be an irreducible affine action 
of the topological group $G$ on ${\H}.$ The linear part of $\alpha$ is trivial  on the FC-center ${\rm FC}(G)$ of
$G;$ more precisely,
every $\ga\in {\rm FC} (G)$ acts as a translation in the direction of ${\cal H}^{\pi(C_\ga)}.$\end{Prop}

\bpr  Let $\ga\in {\rm FC}(G).$ Since $C_\ga$ is a closed subgroup with finite index, by Theorem~\ref{neretin}, 
the restriction of $\alpha$ to $C_\ga$ is irreducible. Hence, by Proposition \ref{Schur},
$\alpha(\ga)$ is a translation by a vector in ${\cal H}^{\pi(C_\ga)}.$ \epr

\medskip
A group $G$ is called an {\it FC-group} if $G= {\rm FC}(G).$
The following result is an immediate consequence of Proposition \ref{FC-center}.

\begin{Prop}\label{FC-irr}
Let $G$ be  an FC-group.
Every irreducible affine action of $G$ on ${\cal H}$, is given by a homomorphism $b:G\rightarrow{\cal H}$ such that ${\rm span}(b(G))$ is dense. \epr
\end{Prop}

We now show that  a result similar to Corollary~\ref{Serre} holds for discrete groups satisfying the following  property introduced
in \cite{LuZ}.

\begin{Def}  A discrete group  $\Ga$ has property (FAb) if, for every  subgroup $H$ of finite index of $\Ga$, we have $\textnormal{Hom}(H,\mathbf{R})=0$.
\end{Def} 

 It is shown in \cite[Proposition 1.30]{LuZ} that $\Ga$ has property (FAb) if and only if $H^1(\Ga,\pi)=0$ for every  complex representation $\pi$ of $\Ga$ with finite image.

\begin{Cor}\label{G/FC(G)} Let $\Ga$ be a group with property (FAb). Then every irreducible affine action $\alpha$ of $\Ga$ factors through $\Ga/{\rm FC}(\Ga)$.
\end{Cor}

\bpr The proof is similar to the proof of Corollary~\ref{G/Z}. \epr

\medskip
We obtain from the previous result the following extension of Serre's result from Corollary~\ref{Serre}, with a similar proof.

\begin{Cor}\label{Serre-FC} Let $\Ga$ be countable discrete group with property (FAb).  If $\Ga/{\rm FC}(\Ga)$ has property (T), then so does $\Ga$. \epr
\end{Cor}

\subsection{Abelian groups}

In this section, $A$ will denote a topological abelian group, written additively. 
Since $A$ is an FC-group, we have from Proposition \ref{FC-irr}, that every
 irreducible affine action of $A$ on $\H$ is given by a continuous homomorphism
 $b:A\to \H$  such that ${\rm span}(b(A))$ is dense. 

\begin{Def} \rm (see \cite{Forst}) A continuous function $Q:A\rightarrow\mathbf{R}^+$ is a {\it non-negative quadratic form} if $Q(x+y)+Q(x-y)=2(Q(x)+Q(y))$ for every $x,y\in A$.
\end{Def}


\begin{Lem}\label{quadratic} A continuous, non-negative function $Q$ on $A$ is a quadratic form if and only if there exists a Hilbert space ${\cal K}$ and a continuous homomorphism $\beta:A\rightarrow{\cal K}$ such that $Q(x)=\|\beta(x)\|^2$ for every $x\in A$.
\end{Lem}

\bpr It is immediate to check that, if $Q(x)=\|\beta(x)\|^2$, then $Q$ is a quadratic form. Conversely, start from a quadratic form $Q$, and observe that $Q(x)=Q(-x)$ and $Q(nx)=n^2Q(x)$ for $n\in\mathbf{N}$ (the latter equality being proved by induction over $n$). Set 
$$V:=A\otimes_\mathbf{Z}\mathbf{R}\qquad \text{and}\qquad \tilde{Q}(x\otimes\lambda)=\lambda^2Q(x);$$ then $\tilde{Q}$ is a well-defined non-negative quadratic form on the real vector space $V$, so we may define ${\cal K}$ as the separation-completion of $V$ and 
$$\beta:A\rightarrow V, \, x\mapsto x\otimes 1$$ does the job. Since the topology of ${\cal K}$ is determined by $Q$ which is continuous, the homomorphism $\beta$ is continuous by construction.\epr

\medskip
Recall that, for an affine isometric action $\alpha$ with cocycle $b$, we denote $\psi(.)=\|b(.)\|^2$.

\begin{Prop}\label{abelian} Let $\alpha$ be an affine action of $A$, with $b(A)$ total in ${\cal H}$. The following properties are equivalent:
\begin{enumerate}
\item[i)] $\alpha$ is irreducible;
\item[ii)] $\psi$ is a quadratic form.
\end{enumerate}
\end{Prop} 

\bpr $(i)\Rightarrow (ii)$ follows immediately from Proposition \ref{FC-irr} and lemma \ref{quadratic}. For $(ii)\Rightarrow (i)$, write $\psi(x)=\|\beta(x)\|^2$, with $\beta:A\rightarrow{\cal K}$ a continuous homomorphism, as in Lemma \ref{quadratic}. Clearly we may assume that $\beta(A)$ is total in ${\cal H}$. The actions $\alpha$ and $\beta$ (viewed as an action by translations) both have total cocycle and define the same function conditionally of negative type, so they are conjugate by an $A$-equivariant affine isometry (see Proposition 2.10.2 in \cite{BHV}).\epr 

\begin{Rem} 
\rm When $A$ is locally compact abelian, it is possible to give a proof of the implication $(i)\Rightarrow(ii)$ in Proposition \ref{abelian}, not depending on Proposition \ref{FC-irr} (so that, together with Lemma \ref{quadratic}, we get a direct proof of Proposition \ref{FC-irr} in the case of an abelian group). Indeed, by the Levy-Khintchine formula (see Theorem 8 in \cite{Forst}), $\psi$ can be written as:
$$\psi(x) = Q(x)+\int_{\hat{A}\backslash\{1_A\}}(1-\mbox{Re}\chi(x))\,d\mu(\chi)$$
where $Q$ is a quadratic form, $\hat{A}$ is the Pontryagin dual of $A$, and $\mu$ is a non-negative measure on $\hat{A}\backslash\{1_A\}$ that gives finite measure to the complement of any neighborhood of the unit $1_A$ of $\hat{A}$. If $\psi$ is not a quadratic form, then $\mu\neq 0$. In this case, choose a point $\chi$ in the support of $\mu$ and a neighborhood $V$ of $\chi$ which is disjoint from some neighborhood of $1_A$. Set then 
$$\psi_0(x)=\int_V(1-\mbox{Re}\chi(x))\,d\mu(\chi), \quad\psi_1(x)=Q(x)+\int_{\hat{A}\backslash(\{1_A\}\cup V)} (1-\mbox{Re}\chi(x))\,d\mu(\chi).$$ Then $\psi=\psi_0+\psi_1$, the functions $\psi_0,\psi_1$ are conditionally of negative type, $\psi_0$ is bounded, and $\psi_0\neq 0$ (because $\mu(V)>0$). By condition (A4) in Proposition \ref{affine}, the action $\alpha$ is reducible.
\end{Rem}

\subsection{Nilpotent groups and FC-nilpotent groups}

The following result generalizes Corollary 5 in \cite{Gui2}, stating that for a nilpotent locally compact group, any non-trivial unitary irreducible representation has zero 1-cohomology.

\begin{Cor}\label{nilp} Let $G$ be a nilpotent group. Any irreducible affine action $\alpha$ of $G$ on ${\cal H}$ is given by a continuous homomorphism $b:G\rightarrow{\cal H}$ such that ${\rm span}(b(G))$ is dense.
\end{Cor}

\bpr We proceed by induction on the nilpotency rank $r$ of $G$, the case $r=1$ being Proposition \ref{FC-irr}. For the general case, let $\alpha$ be an irreducible affine action of $G$, it is enough to show that $\pi$ is the trivial representation, i.e. ${\cal H}^{\pi(G)}={\cal H}$. Assume it is not the case, and let $\alpha_0$ be the projected action on the orthogonal complement of ${\cal H}^{\pi(G)}$. By condition (A6) in Proposition \ref{affine}, the action $\alpha_0$ is irreducible. Since its linear part $\pi_0$ has no non-zero fixed vector, by Proposition \ref{center} the center $Z(G)$ acts trivially in $\alpha_0$, i.e. $\alpha_0$ factors through $G/Z(G)$. By induction hypothesis $\alpha_0$ is an action by translations, meaning that $\pi_0$ is the trivial representation of $G/Z(G)$. This contradiction ends the proof. \epr

\medskip
Denote by $Q$ the convex cone of functions on $G$ of the form $x\mapsto\|b(x)\|^2$, where $b$ is a continuous homomorphism from $G$ to the additive group of a Hilbert space (for $G$ abelian, this is the cone of quadratic forms). 

\begin{Cor} Let $G$ be a nilpotent group. Then $Q$ is the unique maximal face shared by $C$ and $C_u$. \epr
\end{Cor}

\medskip
The  ascending FC-central series $(G_i)_i$ of a group $G$ is defined  inductively
as follows: $G_1={\rm FC}(G)$ and $G_{i+1}$ is the inverse image of ${\rm FC}(G/G_i)$
under the canonical map $G\to G/G_i$  for  every $i\geq 1.$
If  $G_n=G$ and $G_{n-1}\neq G$, then $G$ is said to be \emph{FC-nilpotent}
of rank $n$.
 Examples of FC-nilpotent groups include nilpotent-by-finite groups
and (arbitrary) direct sums   of finite groups.

Corollary~\ref{nilp} cannot be extended to the class of FC-nilpotent
groups. Indeed, $G$ be the semi-direct product $\mathbf{Z}\rtimes C_2$, where
the cyclic group $C_2=\{\pm 1\}$ of order 2 acts on $\mathbf{Z}$ in the non trivial way.
The group $G$ is FC-nilpotent of rank 2; the affine action $\alpha$ of $G$ on $\mathbf{C},$ defined by 
$\alpha(-1, m) x= -x +m$ for $m\in \mathbf{Z}, x\in \mathbf{C},$
is clearly irreducible and not given by a homomorphism $G\to \mathbf{C}.$
Observe that the linear part of $\alpha$ factors though the finite quotient $C_2.$
The next proposition is the proper generalization of this fact.

\begin{Cor}\label{FC-nil} Let $G$ be an  FC-nilpotent and  $\alpha$  an irreducible affine action of $G$ on 
a Hilbert space ${\cal H},$ with linear part $\pi.$ Then $\pi$ can be decomposed
as a direct sum $\pi= \bigoplus_i \pi,$ where each $\pi_i$ is a unitary representation of $G$
which factors through a finite quotient of $G.$
\end{Cor}

\bpr We proceed by induction on the FC-nilpotency rank $r$ of $G$.
When  $r=1,$ the group $G$ is an FC-group and the claim
follows from  Proposition~\ref{FC-irr}. 

Let $r\geq 2.$ Denote by  ${\cal K}$ be the closed linear space of ${\cal H}$ generated
by all subrepresentations of $\pi$ which factor through a finite quotient.
It is clear that the restriction of $\pi$ to ${\cal K}$ can be decomposed as
a direct sum $ \bigoplus_i \pi,$ where each $\pi_i$ is a subrepresentation of $\pi$
which factors through a finite quotient of $G.$

The claim will be proved if we can show that ${\cal K}={\cal H}.$
 Assume, by contradiction, that  this is not the case.
  Let $\alpha_0$ be the projected action on the orthogonal complement ${\cal H}_0$ of ${\cal K}$. 
  By condition (A6) in Proposition\ref{affine}, the action $\alpha_0$ is irreducible. 
 Denote by $\pi_0$ the subrepresentation of $\pi$ defined by ${\cal H}_0.$
  Observe that $\pi_0$   does not factor through a finite quotient of $G.$
  
  Let $\ga\in {\rm FC}(G).$  By Proposition~\ref{FC-center}, 
  $\alpha_0(\ga)$ is a translation in the direction of ${\cal H}_0^{\pi(C_\ga)}.$
  Let   $N_\ga$ be a  normal subgroup of finite index of $G$  contained in $C_\ga.$
  Then $ {\cal H}_0^{\pi(N_\ga)}$  is a $\pi(G)$-invariant 
  subspace of ${\cal H}_0$ and the corresponding subrepresentation
  of $\pi_0$ factors through the finite quotient $G/N_\ga.$
  It follows that $ {\cal H}_0^{\pi(N_\ga)}=\{0\}$ 
  and hence ${\cal H}_0^{\pi(C_\ga)} =\{0\}.$ 
  So, $\alpha_0(\ga)$ is the identity.
  We have therefore proved that $\alpha_0$  factors
   through $G/{\rm FC}(G)$. 
   
   Observe that $G/{\rm FC}(G)$ is FC-nilpotent
   of rank $ r-1.$
   By induction hypothesis,  $\pi_0$ is a direct sum of subrepresentations which factor though finite quotients;
   hence, ${\cal H}_0=\{0\}$ and  this is a contradiction. \epr

\subsection{The left regular representation of  a discrete group}
\label{SS: LeftRegRep}
For a discrete group $\Ga,$ 
we will be interested in the  question of the existence of an irreducible affine isometric action with linear part 
the left regular representation $\la_\Ga.$
More generally, we will consider the same question for  a closed $\Ga$-invariant subspace $\H$ of a countably many copies of $\ell^2(\Ga);$
thus, $\H$ is  a closed subspace of $\oplus_{n\in \NN} \ell^2(\Ga)$ which
is invariant under the representation $\oplus_{n\in \NN}\la_\Ga.$ 
Observe that such a space $\H$ is a Hilbert module over the  left group von Neumann algebra $L(\Ga)$
and every  Hilbert module over $L(\Ga)$ is of this form (see below).

Let $\M$ be finite von Neumann algebra, that is, $\M$ is a von Neumann algebra equipped with 
a faithful normal  finite trace $\tau: \M\to \mathbf C$.
Let $L^2(\M)$ be the Hilbert space  obtained from $\tau$ by 
the GNS construction. We  identify $\M$  with the subalgebra of ${\B}(L^2({\M}))$ of operators
given by left multiplication with elements from $\M.$ The commutant of $\M$ in
$\B(L^2(\M))$  is $\M '=J\M J,$ where $J : L^2(\M) \to L^2(\M)$ is the conjugate linear isometry
which extends the mapping $\M\to \M, x\mapsto x^*$. The trace on  $\M'$, again denoted by $\tau,$
is defined by  $JxJ\mapsto \tau (x)$ for $x\in \M.$
	
Let $\H$ be a  Hilbert $\M$-module, that is, a separable Hilbert space  with a unital normal homomorphism 
$\M\to {\B}({\H}).$ Then   $\H$ can be identified as $\M$-module  to a submodule of
$L^2(\M)\otimes \K$ for an infinite dimensional separable Hilbert space
$\K,$ where $\M$ acts on  $L^2(\M)\otimes {\K}$  by  
$$
\xi\otimes \eta \mapsto T\xi \otimes \eta, \qquad T\in{\M}, \xi\in L^2(\M), \eta \in {\K}.
$$

Let $P: L^2({\M})\otimes {\K} \to {\H}$ be the orthogonal projection.
Then $P$ belongs to the commutant of $\M$
in ${\B}(L^2({\M})\otimes {\K})$, which is  ${\M}' \otimes {\B}({\K})$, where
${\M}'$ is the commutant of $\M$ in ${\B}(L^2({\M})).$

Let $\{e_n\}_n$ be a Hilbert space basis of $\K.$ Let $(P_{ij})_{i,j}$ be the matrix of 
$P$ with respect to the decomposition  $L^2({\M})\otimes {\K} = \oplus_{i} (L^2({\M}) \otimes \CCC e_i).$
Then each $P_{ij}$ belongs to ${\M} '.$
The von Neumann dimension  of the $\M$-module $\H,$ which takes values in $[0, +\infty[ \cup \{+\infty\}$,
is defined by
$$
 {\rm dim}_{\M}{\H} = \sum_{i} \tau(P_{ii}).
  $$
When $\M$ is a factor, $\H$ is characterized as $\M$-module by its von Neumann dimension,
up to unitary equivalence (see e.g. Proposition~3.2.5 in \cite{GHJ}).

Let $\Ga$ be a discrete countable group and $\la_\Ga$  the left regular representation of $\Ga$ on $\ell^2(\Ga).$
Denote by $L(\Ga)$ the left regular von Neumann algebra of $\Ga.$
Recall that $L(\Ga)$ is the closure 
of the linear span of  $\{\la_\Ga(\ga)\ : \ \ga \in \Ga\}$ in the weak (or strong) operator topology.
The commutant  $L(\Ga)'$  of $L(\Ga)$ in ${\B}(\ell^2(\Ga))$ is the right group von Neumann algebra $R(\Gamma),$ the von Neumann algebra
generated by the right regular representation of $\Gamma$.
The algebras $L(\Ga)$ and $R(\Ga)$ are finite
von Neumann algebras:  a faithful normal trace $\tau$ on $L(\Gamma)$ or $R(\Ga)$ is given by 
$$
\tau(T) =\langle T\delta_e| \delta_e\rangle, 
 \qquad \mbox{for all}\quad T\in L(\Ga)\quad\mbox{or}\quad T\in R(\Ga).
 $$

Assume now that $\Ga$ is non amenable and finitely generated.  By  \cite{BV},  there exists a $R(\Ga)$-equivariant isomorphism between the first cohomology $H^1(\Gamma,\lambda_\Gamma)$ and the first $L^2$-cohomology $H^1_{(2)}(\Gamma)$; it follows that
 $H^1(\Ga; \la_\Ga)$  has a Hilbert space structure.
The first $L^2$-Betti number of $\Ga$  
is 
$$\beta^1_{(2)}(\Gamma) = \dim_{R(\Gamma)}H^1_{(2)}(\Gamma).$$

Recall that 
$L(\Ga)$ or $R(\Ga)$  is a factor (that is, their common center consists only of the scalar multiples of the identity) if and only if $\Ga$ is ICC, i.e. every non-trivial conjugacy class in $\Ga$ is infinite; otherwise said, ${\rm FC}(\Ga)$ is trivial.

The following result was initially obtained in the special case of the $L(\Ga)$-module $\ell^2(\Ga)$ under the 
additional assumption that $\Ga$ is an ICC group; we thank S. Vaes for suggesting to jack it up to arbitrary $L(\Ga)$-modules.

\begin{Thm}
\label{AffRegRep}
 Let $\Gamma$ be a non-amenable, finitely generated group,  and let
$\H$ be 	a non zero Hilbert $L(\Ga)$-module with finite von Neumann dimension.
 Denote by 
 $\la^{\H}$ the corresponding  unitary representation of   $\Ga$ in  $\H.$
The following properties are equivalent:
\begin{enumerate}
\item[i)] there exists an irreducible affine isometric action of $\Ga$ with linear part $\la^{\H}$;
\item[ii)] ${\rm FC}(\Ga)$ is finite, ${\rm FC}(\Ga)$ acts trivially on $\H,$ 
and  
$$\beta^1_{(2)} (\Ga/ {\rm FC}(\Ga)) \geq {\rm dim}_{L(\Ga/{\rm FC}(\Ga))} \H.$$  
\end{enumerate}
\end{Thm}

\bpr
\mbox{\bf First step:} we assume that $\Ga$ is an ICC group, so that 
$L(\Ga)$ is a factor. 

Since ${\rm dim}_{L(\Ga)}{\H}$ is finite,  we can find an integer $k$  such that $\H$ is  a submodule of
$\ell^2(\Ga)\otimes \mathbf{C}^k$.

 Let $P: \ell^2(\Ga)\otimes \mathbf{C}^k\to \H$
be the corresponding orthogonal projection  with range $\H.$ 
Set  ${\M}= L(\Ga) \otimes I_{\mathbf{C}^k} \cong L(\Ga)$.
The commutant  of ${\M}$ in ${\B}(\ell^2(\Ga)\otimes \mathbf{C}^k)$ is 
$${\M}' =R(\Ga) \otimes {\B}(\mathbf{C}^k) =M_k(R(\Ga)).$$
So, we  can write 
$$P= (P_{ij})_{1\leq i,j\leq k} \in R(\Ga) \otimes {\B}(\mathbf{C}^k) =M_k(R(\Ga))
$$  
and 
$$
 {\rm dim}_{L(\Ga)}{ \H}=\sum_{i=1}^k \tau(P_{ii})
  $$
The subalgebras ${\M} P$ and $P{\M}' P$ of  $\B(\H) $ are finite factors
and  we have  $P {\M}' P = ({\M} P)';$ thus, the commutant 
of $\la^{\H}(\Ga)$ is $P {\M}' P.$

Next,  since $\Ga$ is not amenable, the $1$-cohomology group 
$H^1(\Ga, \oplus_{i=1}^k\la_\Ga)$ coincides with the \emph{reduced cohomology group} 
$\overline{H}^1(\Ga, \oplus_{i=1}^k\la_\Ga)$, that is, the quotient of $Z^1$ by the closure of $B^1$,
for the topology of pointwise convergence of $\Ga$ (\cite[Corollaire1]{Gui2}); moreover, we have
$$
\overline{H}^1(\Ga, \oplus_{i=1}^k \la_\Ga)= \oplus_{i=1}^k\overline{H}^1(\Ga, \la_\Ga)= H^1_{(2)}(\Gamma)\otimes\mathbf{C}^k,
$$
which is a module over $\M'.$
It follows that the 1-cohomology of $\la^{\H}$ 
 is  given by the $P \M' P $-module  $P (H^1_{(2)}(\Gamma)\otimes\mathbf{C}^k).$

By Corollary \ref{separating},  there exists an irreducible affine isometric action of $\Ga$ with linear part $\la^{\H}$
if an only if  $P (H^1_{(2)}(\Gamma)\otimes\mathbf{C}^k)$ admits a separating vector for 
$P \M' P.$
Now,  ${\rm dim}_{P{\M}' P} P (H^1_{(2)}(\Gamma)\otimes\mathbf{C}^k)$ is the   coupling
constant for $P{\M}' P$ acting on $P (H^1_{(2)}(\Gamma)\otimes\mathbf{C}^k);$ 
see \cite[Proposition 3.2.5]{GHJ}.
Hence
$P (H^1_{(2)}(\Gamma)\otimes\mathbf{C}^k)$ admits a separating vector for 
$P {\M}' P $ if only if  
$$
{\rm dim}_{P{\M}' P} P (H^1_{(2)}(\Gamma)\otimes\mathbf{C}^k) \geq 1
$$
(see    \cite[Chap. III, \S 6, Proposition 3]{DixmierVN}).

On the other hand, by \cite[Chap. III, \S 6, Proposition 2]{DixmierVN} or \cite[Proposition 3.2.5]{GHJ}, we
have 
$$
{\rm dim}_{P{\M}' P} P (H^1_{(2)}(\Gamma)\otimes\mathbf{C}^k)\delta_{\M'}(P)
={\rm dim}_{\M'}(H^1_{(2)}(\Gamma)\otimes\mathbf{C}^k),$$
where $\delta_{\M'}$ is the canonical normalized trace on ${\M}'=M_k(R(\Ga))$.
We have, for every $T= (T_{ij})_{1\leq i,j\leq k} \in M_k(R(\Ga)),$
$$
 \delta_{\M'}(T)= \frac{1}{k} \sum_{i=1}^k \tau(T_{ii})
 $$
 and hence
  $$
 \delta_{\M'}(P)= \frac{1}{k}  {\rm dim}_{L(\Ga)} {\H}.
 $$
 Moreover 
$$
{\rm dim}_{M_k(R(\Ga))}(H^1_{(2)}(\Gamma)\otimes\mathbf{C}^k)=\frac{{\rm dim}_{R(\Ga)}H^1_{(2)}(\Gamma)}{k}=\frac{\beta^1_{(2)}(\Gamma)}{k}.
$$
We have therefore
$$
{\rm dim}_{P{\M}' P} P (H^1_{(2)}(\Gamma)\otimes\mathbf{C}^k)  {\rm dim}_{L(\Ga)} {\H} = 
\beta^1_{(2)}(\Gamma).
$$
As a consequence, 
$${\rm dim}_{P{\M}' P} P (H^1_{(2)}(\Gamma)\otimes\mathbf{C}^k)\geq 1$$
 if and only
if   $\beta^1_{(2)}(\Gamma) \geq  {\rm dim}_{L(\Ga)} \H.$
 
\noindent
 \mbox{\bf Second step:} we  assume that ${\rm FC}(\Ga)$ is non trivial.
 Observe that $\Ga/{\rm FC}(\Ga)$ is not amenable,
 since ${\rm FC}(\Ga)$ is amenable and $\Ga$ is not amenable.
 
  Assume first that  there exists an irreducible affine isometric  action $\alpha$ of $\Ga$ 
 with linear part $\la^{\H}$.
By Proposition~\ref{FC-center}, $\la^{\H}$ is trivial on  ${\rm FC}(\Ga).$ 
Since $\la^{\H}$ is a subrepresentation of  a multiple of the regular representation $\la_\Ga,$
 it follows that  ${\rm FC}(\Ga)$ is finite.  
 As a consequence, 
$\ell^2(\Ga/{\rm FC}(\Ga))$ can be identified as  $L(\Ga)$-module (or as $R(\Ga)$-module)
with the closed subspace  $\ell^2(\Ga)^{\la_\Ga ({\rm FC}(\Ga))}$  of $\ell^2(\Ga).$  So, the
Hilbert module  $\H$ over $L(\Ga)$, on which ${\rm FC}(\Ga)$ acts trivially,  can be
identified with a Hilbert module over  $L(\Ga/{\rm FC}(\Ga)).$

Since ${\rm FC}(\Ga)$ is finite, it is straightforward to check that $\Ga/{\rm FC}(\Ga)$ is  ICC.
By the first step, it follows that 
$$\beta^1_{(2)} (\Ga/{\rm FC}(\Ga))\geq {\rm dim}_{L(\Ga/{\rm FC}(\Ga))} \H.$$

Conversely, assume that ${\rm FC}(\Ga)$ is finite, that ${\rm FC}(\Ga)$ acts trivially on $\H,$ 
and  that
$$\beta^1_{(2)} (\Ga/{\rm FC}(\Ga)))\geq {\rm dim}_{L(\Ga/{\rm FC}(\Ga))} \H.$$
It follows by the first step that  there exists an irreducible affine isometric action of $\Ga/{\rm FC}(\Ga)$ with linear part 
given by $\la^{\H}.$
This concludes the proof.\epr

\medskip

As a corollary, we obtain a necessary condition for the existence
 an irreducible affine isometric action of $\Ga$ with linear part $\la^{\H}$,
 in terms  of $\beta^1_{(2)} (\Ga)$ and ${\rm dim}_{L(\Ga)} \H.$
\begin{Cor}\label{Cor1-AffRegRep} 
Let   $\Gamma$, $\H$ and $\la^{\H}$ be as in Theorem~\ref{AffRegRep}.
If there  exists an irreducible affine isometric action of $\Ga$ with linear part $\la^{\H}$,
then 
$$\beta^1_{(2)} (\Ga)\geq {\rm dim}_{L(\Ga)} \H.$$
\end{Cor}
\bpr
 By Theorem~\ref{AffRegRep}, the cardinality $N$ of ${\rm FC}(\Ga)$ is finite. It is easily checked that 
${\rm dim}_{L(\Ga/{\rm FC}(\Ga))} {\H}= N {\rm dim}_{L(\Ga)} {\H}$; similarly, 
since $H^1_{(2)}(\Ga/{\rm FC}(\Ga))$ can be identified with the $R(\Ga)$-submodule
of $H^1_{(2)}(\Ga)$ on which ${\rm FC}(\Ga)$ acts trivially, we have 
$$N\beta^1_{(2)} (\Ga) \geq \beta^1_{(2)} (\Ga/{\rm FC}(\Ga))$$
and hence, using   Theorem~\ref{AffRegRep}, we obtain 
$$\beta^1_{(2)} (\Ga)\geq {\rm dim}_{L(\Ga)} \H. 
$$
 \epr
 
 The following corollaries are immediate consequences of  Theorem~\ref{AffRegRep}.
\begin{Cor}\label{No-ICC} Let $\Gamma$ be a non-amenable, finitely generated group
such that ${\rm FC}(\Ga)$ is infinite. No non-zero
 $L(\Ga)$-module  $\H$ has an irreducible affine isometric action with linear part $\la^{\H}.$
\end{Cor}

\begin{Cor}\label{ICCnewbetti} For $\Gamma$ a non-amenable, finitely generated ICC group, we have
$$\beta_{(2)}^1(\Ga)=\sup\{t \geq 0: t.\lambda_\Gamma\;\mbox{is the linear part of an irreducible affine action}\},$$
where $t.\lambda_\Ga$ is the underlying $\Ga$-representation of  the unique $L(\Ga)$-module of von Neumann dimension $t.$
\end{Cor}

\begin{Ex} \rm 
(i) The group $PSL_2(\mathbf{Z})$ is ICC and satisfies $\beta_{(2)}^1(PSL_2(\mathbf{Z}))=\frac{1}{6}$ (see Section 4 in \cite{CG}), so there exists no irreducible affine action with linear part the left regular representation.

\noindent 
(ii) Let $\widetilde G$ be the universal cover of the Lie group
$G= SL_2(\mathbf{R})$ and let $\Ga$ be the inverse image in $\widetilde G$
of $SL_2(\mathbf{Z})$
under the covering map  $\widetilde G\to G.$ Then, since ${\rm FC}(\Ga)$ is infinite, 
no non-zero $L(\Ga)$-module  $\H$ has an irreducible affine isometric action with linear $\la^{\H}.$
\end{Ex}

For the free group $\mathbf{F}_n$ on $n$ generators ($2\leq n\leq+\infty$), we have $\beta_{(2)}^1(\mathbf{F}_n)=n-1$ (see  \cite{CG})
and it is possible to construct explicit irreducible affine isometric actions with linear part $\lambda_{\mathbf{F}_n}$. Indeed, let $(a_i)_{1\leq i\leq n}$ be a free generating family of $\mathbf{F}_n$. Set $b(a_1)=\delta_1$ (the characteristic function of the identity of $\mathbf{F}_n$), and $b(a_i)=0$ for $i\geq 2$. Since $\mathbf{F}_n$ is free, we may extend uniquely $b$ to a 1-cocycle $b\in Z^1(\mathbf{F}_n, \lambda_{\mathbf{F}_n})$. It is easily seen that, for $k\geq 0$, we have $b(a_1^k)=\sum_{i=0}^{k-1}\delta_{a_1^i}$, so that $b$ is unbounded.

\begin{Prop} For $b$ as above, the affine isometric action of $\mathbf{F}_n$ on $\ell^2(\mathbf{F}_n)$ given by $\alpha(g)v=\lambda_{\mathbf{F}_n}(g)v+b(g)$, is irreducible.
\end{Prop}

\bpr Let $Av=Tv+t$ be an affine transformation of $\ell^2(\mathbf{F}_n)$ in the commutant of $\alpha$. Then $T\in R(\mathbf{F}_n)$ and $(T-1)b(g)=\lambda_{\mathbf{F}_n}(g)t-t$ for every $g\in\mathbf{F}_n$. For $g=a_2$, we get $\lambda_{\mathbf{F}_n}(a_2)t=t$, hence $t=0$ since $a_2$ has infinite order. So $(T-1)b(g)=0$ for every $g$. For $g=a_1$, this gives $(T-1)\delta_1=0$, hence $T=1$ since $\delta_1$ is separating for $R(\mathbf{F}_n)$. By Proposition \ref{Schur}, the action $\alpha$ is irreducible. \epr

\medskip
The situation is completely different for the regular representation of amenable groups. Indeed we have the following result due to Andreas Thom, who kindly gave us permission to include it here.

\begin{Thm}\label{thom} Let $\Gamma$ be a discrete, amenable group. Let $\alpha$ be an affine isometric action of $\Gamma$, with linear part $\lambda_\Gamma$. For every $\varepsilon>0$, the action $\alpha$ admits a closed, affine invariant subspace ${\cal H}_\varepsilon$ such that the linear part ${\cal H}_\varepsilon^0$ satisfies $\dim_{L(\Gamma)} {\cal H}^0_\varepsilon<\varepsilon$. In particular, there is no irreducible affine action of $\Gamma$ with linear part $\lambda_\Gamma$.
\end{Thm}

Observe that, by a result of Guichardet \cite[Corollaire1]{Gui2}, we have $H^1(\Gamma,\lambda_\Gamma)\neq 0$ for every countable amenable group $\Gamma$, so there is indeed something to be proved.

\medskip
\bpr Let $b\in Z^1(\Gamma,\lambda_\Gamma)$ be the 1-cocycle defining $\alpha$. We will need the ring ${\cal U}(\Gamma)$ of operators affiliated to the von Neumann algebra $R(\Gamma)=\lambda_\Gamma(\Gamma)'$, as introduced e.g. in \cite[Chap.8]{Lue}. We recall that, as $\Gamma$-modules, we have the chain of inclusions $R(\Gamma)\subset\ell^2(\Gamma)\subset{\cal U}(\Gamma)$. Now we appeal to a special case of Theorem 2.2 in \cite{PeTh}: if a group $\Lambda$ has vanishing first $L^2$-Betti number, then $H^1(\Lambda,{\cal U}(\Lambda))=0$. This applies to $\Gamma$, by the Cheeger-Gromov vanishing theorem for amenable groups (Theorem 0.2 in \cite{CG}). This means that, viewing our cocycle $b\in Z^1(\Gamma,\ell^2(\Gamma))$ as a cocycle in $Z^1(\Gamma,{\cal U}(\Gamma))$, we may trivialize it and find some $f\in{\cal U}(G)$ such that $b(g)=\lambda_\Gamma(g)f-f$ for every $g\in\Gamma$. We now proceed as is the proof of Corollary 2.4 in \cite{PeTh}: given $\varepsilon>0$, we find a projector $Q\in R(\Gamma)$ such that $Qf\in\ell^2(\Gamma)$ and $\dim_{R(\Gamma)} (1-Q)(\ell^2(\Gamma))<\varepsilon$. It is then easy to check (as in the proof of our Proposition \ref{affine}) that the closed affine subspace ${\cal H}_\varepsilon=:-Qf + (1-Q)(\ell^2(\Gamma))$ is $\alpha(\Gamma)$-invariant.\epr

\section{Direct sums of irreducible actions}

For affine isometric actions $\alpha_1,\alpha_2$ of a group $G$, we may consider in an obvious way the direct sum $\alpha_1\oplus\alpha_2$. Unlike the direct sum of unitary representations, which is always reducible, it may happen that the direct sum of two affine isometric actions is irreducible. For instance, if $\beta_1,\beta_2$ are linearly independent homomorphisms $G\rightarrow\mathbf{C}$, then $\beta_1\oplus\beta_2$ defines an irreducible affine isometric action of $G$ on $\mathbf{C}^2$. On the other hand, if $\alpha$ is any affine isometric action of $G$, then $\alpha\oplus\alpha$ is not irreducible (look at the diagonal). We shall give a sufficient and necessary condition for the direct sum of two irreducible actions to be irreducible. 

\medskip
In order to state the main result of this section (Theorem \ref{DirectSums} below) we need to clarify the notion of equivalence between affine isometric actions.
\begin{Def} \rm
Let $\alpha_1$ and $\alpha_2$ be two affine isometric actions of a group $G$. We say that $\alpha_1$ and $\alpha_2$ are equivalent if they are intertwined by an invertible continuous affine mapping, that is, if there exists an invertible continuous affine mapping $A:\H_{\alpha_1}\rightarrow\H_{\alpha_2}$ satisfying:
$$A\alpha_1(g)=\alpha_2(g)A, \quad\text{for all} \, g\in G.$$
\end{Def}

If we write $A(\cdot)=T(\cdot)+t$ and $\alpha_i(g)(\cdot)=\pi_i(g)(\cdot)+b_i(g)$, the above definition boils down to $T\pi_1(g)=\pi_2(g)T$ and $Tb_1(g)=b_2(g)+\pi_2(g)t-t$ for all $g\in G.$

\medskip
Since the actions are by isometries, it may seem more natural to require the intertwining in the definition of equivalence to be given by an isometric operator, in which case we would say that  the actions are  isometrically equivalent.
To motivate our definition, one should be reminded of  the similar definition for unitary representations. It is well-known that, in this case, an equivalence can always be implemented via a unitary intertwiner. This is a consequence of the fact that every invertible intertwiner can be ``straightened" by replacing it with its unitary part (see e.g. \cite[Appendix A.1]{BHV}). However, this fails for affine isometric actions:
equivalent affine actions by isometries need not be  isometrically equivalent
\footnote{As an example, consider two actions of $\mathbf{Z}$ on $\mathbf{R}$, the first one by integer translations, the second one by even translations. These actions are equivalent in our sense, but clearly they are not isometrically equivalent.}.

\begin{Thm}\label{DirectSums} Let $\alpha_1,\alpha_2$ be irreducible affine isometric actions of a group $G$. The following properties are equivalent:
\begin{enumerate}
\item[i)] $\alpha_1\oplus\alpha_2$ is reducible.
\item[ii)] $\alpha_1$ and $\alpha_2$ admit equivalent projected actions.
\end{enumerate}
\end{Thm}

Before proving this theorem, we pinpoint two specific cases, important enough to be considered on their own.

\medskip
Recall that two unitary representations $\pi,\sigma$ of $G$ are said to be {\it disjoint} if ${\rm Hom}_G({\cal H}_\pi,{\cal H}_\sigma)=0$.

\begin{Prop}\label{directsum} Let $\alpha_1,\dots,\alpha_k$ be irreducible affine actions of $G$, with linear parts $\pi_1,\dots,\pi_k$. Assume that the $\pi_i$'s are pairwise disjoint. Then the direct sum $\alpha:=\alpha_1\oplus\dots\oplus\alpha_k$ is irreducible. 
\end{Prop}

\bpr Let $b=(b_1,\dots,b_k)$ be the 1-cocycle defining $\alpha$. Let $Av=Tv+t$ be a continuous affine mapping in the commutant of $\alpha$. Write $T$ as a $k\times k$-matrix $(T_{ij})_{1\leq i,j\leq k}$ where $T_{ij}$ is a bounded operator ${\cal H}_{\pi_j}\rightarrow{\cal H}_{\pi_i}$; similarly, write $t=(t_1,\dots,t_k)$. Since $T$ belongs to the commutant of $\pi_1\oplus\dots\oplus\pi_k$, we have $T_{ij}\in {\rm Hom}_G({\cal H}_{\pi_j},{\cal H}_{\pi_i})$
and  hence $T_{ij}=0$ for $i\neq j$. The relation $(T-1)b(g)=\partial_t(g)$ then gives 
$$(T_{ii}-1)b_i(g)=\partial_{t_i}(g) \quad\text{for}\ 1\leq i\leq k \quad \text{and}\ g\in G.$$
 This means that the affine map $A_iv=:T_{ii}v+t_i$ is in the commutant of $\alpha_i$. Since the latter is irreducible, we get $T_{ii}=1$; hence $T=1$ and $\alpha$ is irreducible.
\epr

\medskip
For $\pi$ a unitary representation of $G$ and $k\in\mathbf{N}$, we denote by $k\cdot\pi$ the representation $\pi\oplus\dots\oplus\pi$ ($k$ times). ).

\begin{Prop}\label{linearlyindep} Let $\pi$ be an irreducible unitary  representation of $G$. Let $b_1,\dots,b_k$ be elements in $Z^1(G,\pi)$ whose classes $[b_1],\dots,[b_k]$ are linearly independent in $H^1(G,\pi)$. Then the affine isometric action $\alpha=\bigoplus_{i=1}^k \alpha_{\pi,b_i}$ is irreducible. 
\end{Prop}

\bpr Let $Av=Tv+t$ be a  continuous affine mapping in the commutant of $\alpha$. In view of Proposition \ref{Schur}, we have to show that $A$ is a translation, that is, $T=1$. We know that $T$ is in the commutant of $k\cdot\pi$ and that $(T-1)b=\partial t$, where $b=\oplus_{i=1}^k b_i$.

Write $T$ as a $k\times k$-matrix $(T_{ij})_{1\leq i,j\leq k},$ where $T_{ij}$ is a bounded operator ${\cal H}_{\pi}\rightarrow{\cal H}_{\pi}.$
Then every $T_{ij}$ interwines $\pi$ with itself and hence $T_{ij}=\la_{ij}1$ for some $\la_{ij}\in {\mathbf C},$ by  Schur's lemma.
On the other hand, since $$H^1(G,k\cdot\pi)=\underbrace{H^1(G,\pi)\oplus\dots\oplus H^1(G,\pi)}_{k\, \textrm{times}},$$
 we have
$$(T-1)\left(\begin{array}{c} [b_1] \\  \vdots \\ \mbox{$[b_{k}]$} \end{array}\right)=0;$$
 since the $[b_i]$'s are linearly independent, we deduce that $T=1$.
\epr

\medskip

\noindent \textbf{Proof of Theorem \ref{DirectSums} }: ~  Denote by $\pi_1, b_1$ and $\pi_2, b_2$ the linear and translation parts of the actions
 $\alpha_1$ and $\alpha_2$.
\medskip

$(ii)\Rightarrow (i)$ There exist non zero  $(\pi_1\oplus\pi_2)(G)$-invariant closed linear subspaces ${\cal K}_1$ and ${\cal K}_2$ of $\mathcal\H_{\pi_i}$ such that the projected actions of $\alpha_1$ and $\alpha_2$ on ${\cal K}_1$ and ${\cal K}_2$ are equivalent. Let $A:{\cal K}_1\rightarrow{\cal K}_2$ be a continuous affine, invertible map implementing the equivalence. Then the graph of $A$ is a proper closed, invariant, affine subspace of the projected action of $\alpha_1\oplus\alpha_2$ onto $\mathcal\K_1\oplus\mathcal\K_2$. This contradicts characterization (A6) of irreducibility from Proposition~\ref{affine}.

$(i)\Rightarrow (ii)$ Since $\alpha_1\oplus\alpha_2$ is reducible, we can find, by (A3) from Proposition~\ref{affine},   a non-zero closed linear subspace $\cal K$ of $\H_{\pi_1}\oplus\H_{\pi_2}$ which is invariant under  $(\pi_1\oplus\pi_2)(G)$ and
 such that the projection of $b=b_1\oplus b_2$ on $\K$ is bounded. Upon  conjugating $\alpha=\alpha_1\oplus\alpha_2$ by a translation, we may assume that the projection  of $b$ on $\K$ is $0$. 
 
 Denote by $P_i:\: \K\rightarrow \H_{\pi_i}$ the orthogonal projection of $\K$ onto $\H_{\pi_i}$. We may also assume that 
 $P_i(\K)$ is dense in $\H_{\pi_i}$ for $i=1,2;$ indeed, otherwise we can replace $\alpha$ by its projected action
 on   $\overline{P_1(\K)}\oplus \overline{P_2(\K)}.$ 
 
 Next, observe that $\K$ is transverse to the $\H_{\pi_i}$'s. Indeed, if the intersection $\K\cap\H_{\pi_i}$ were non-zero, the projection of $b_i$ 
 on $\K\cap\H_{\pi_i}$ being bounded, this would contradict the irreducibility of $\alpha_i$. 
So, $P_1$ and $P_2$ are injective.  We can therefore consider the densely defined, unbounded,  invertible closed operator $S=P_2P_1^{-1}$
(for background about unbounded operators, see e.g. \cite[Chap. 5]{Ped2}). 
Note that $\K$ being $(\pi_1\oplus\pi_2)(G)$-invariant, it is immediate that the domain ${\mathcal D}(S)$ of $S$ is $\pi_1(G)$-invariant, that its range is $\pi_2(G)$-invariant and that $S$ intertwines the corresponding two subrepresentations of $\pi_1$
and $\pi_2$ (on non-closed subspaces!). Now, recall that,
for every $g\in G$, the vector $b(g)=b_1(g)\oplus b_2(g)$ is orthogonal to $\K$;  hence, we have
$$\langle b_1(g),v\rangle+\langle b_2(g),Sv\rangle=0\quad \text{for all} \ v\in\mathcal D(S).$$
This relation implies that $$|\langle b_2(g),Sv\rangle|=|\langle b_1(g),v\rangle|\le\|b_1(g)\|\|v\|;$$
 hence $b_2(g)$ belongs to the domain of $S^*$ and $b_1(g)=-S^*b_2(g)$
 for all $g\in G.$  This shows that  $-S^\star$ intertwines $\alpha_2,$ projected on the domain of $S^*,$ and $\alpha_1$.
 
The closed operator  $S^\star$ has a  polar decomposition  $-S^\star=UT$,  where $U:\H_{\pi_2}\rightarrow\H_{\pi_1}$ is unitary and $T:{\mathcal D}(S)\rightarrow\H_{\pi_2}$ is a positive  unbounded closed operator.
Let $B$ be a bounded Borel subset of the spectrum of $T$ with positive measure, and denote by $P_B$ the corresponding spectral projector. 
Then $-S^\star P_B$ is a bounded operator and provides an equivalence between  $\alpha_2$ projected on ${\rm Im}(P_B)$ and $\alpha_1$ projected on ${\rm Im}(S^\star P_B)$. This concludes the proof.\epr

\section{Products and lattices in products}
\subsection{Product groups}

 The following result about irreducible affine actions of product groups
is a consequence of a result of Shalom from \cite{ShaInv} and Proposition \ref{directsum}.
\begin{Prop}\label{irredproduct} Let $G=G_1\times\dots\times G_n$ be the product of non-trivial, compactly generated, locally compact groups. Let $\pi$ be a unitary representation of $G$, not weakly containing the trivial representation, and let $\alpha$ be an affine isometric action of $G$ with linear part $\pi$. The following properties are equivalent:
\begin{enumerate}
\item[i)] $\alpha$ is irreducible.
\item[ii)] $\alpha\simeq\alpha_1\oplus\dots\oplus\alpha_n$, where $\alpha_i$ is an irreducible affine action of $G$ factoring through $G_i$ for every $i=1,\dots\,n$.
\end{enumerate}
\end{Prop}

\bpr $(i)\Rightarrow(ii)$ Set $H_i=G_1\times\dots\times G_{i-1}\times\{1\}\times G_{i+1}\times\dots\times G_n$. Let $b\in Z^1(G,\pi)$ be the cocycle defining $\alpha$. We appeal to a result of Shalom (\cite{ShaInv}, Theorem 3.1; this uses the assumption that $\pi$ does not weakly contain the trivial representation): $b$ is cohomologous to a sum $b_1 + \dots+b_n$, where $b_i$ is a cocycle factoring through $G_i$ and taking values in the space ${\cal H}^{\pi(H_i)}$ of $\pi(H_i)$-fixed vectors. Upon  conjugating $\alpha$ by a translation, we may assume that $b=b_1+\dots+b_n$. 
Denote by $\pi_i$ the subrepresentation of $\pi$ defined by 
the $\pi(G)$-invariant  space  ${\cal H}^{\pi(H_i)}.$ Since $\pi_i$ factors through $G_i,$ the only possible common sub-representation of $\pi_i$ and $\pi_j$ for $i\neq j$ is the trivial representation, which is ruled out by the fact that $\pi$ has no non-zero fixed vector.
Hence, the spaces ${\cal H}^{\pi(H_i)}$ are pairwise orthogonal, so $b=b_1\oplus \dots\oplus b_n$. By irreducibility of $\alpha$, we have ${\cal H}={\cal H}^{\pi(H_1)}\oplus\dots\oplus{\cal H}^{\pi(H_n)}$.  

Define $\alpha_i$ as the projected action of $\alpha$ on ${\cal H}^{\pi(H_i)}$. By construction, $\alpha=\alpha_1\oplus\dots\oplus\alpha_n$ and $\alpha_i$ factors through $G_i$; finally $\alpha_i$ is irreducible, by (A6) from Proposition \ref{affine}.

$(ii)\Rightarrow(i)$ Let $\pi_i$ be the linear part of $\alpha_i$. As above,  the $\pi_i$'s are pairwise disjoint representations of $G$, since $\pi_i$ factors through $G_i$. So Proposition \ref{directsum} applies, and $\alpha$ is irreducible.\epr

\begin{Cor}\label{irredproduct2} Keep notations as in Proposition \ref{irredproduct}. Let $\pi$ be an  irreducible unitary representation of $G$, not weakly containing the trivial representation. If $H^1(G,\pi)\neq 0$, then $\pi$ factors through $G_i$ for some $i=1,\dots, n$.
\end{Cor}

\bpr Let $b\in Z^1(G,\pi)$ be a cocycle which is not a coboundary. By Example \ref{basicex2}, the affine action $\alpha_{\pi,b}$ is irreducible. By Proposition \ref{irredproduct}, we can write $\alpha=\alpha_1\oplus\dots\oplus\alpha_n$, where $\alpha_i$ factors through $G_i$. Let $\pi_i$ be the linear part of $\alpha_i$, so that $\pi=\pi_1\oplus\dots\oplus\pi_n$. By irreducibility of $\pi$, only one of the $\pi_i$'s can be
a non-zero representation.\epr

\medskip 
We note that the assumption that $\pi$ does not weakly contain the trivial representation is necessary in Proposition \ref{irredproduct} and Corollary \ref{irredproduct2}. To see it, let us introduce, 
for a discrete group $\Gamma$, the ``left-right'' representation $\vartheta$ of $\Gamma\times\Gamma$ on $\ell^2(\Gamma)$, defined by:
$$(\vartheta(g,h)\xi)(x)=\xi(g^{-1}xh), \quad  \xi\in\ell^2(\Gamma), g,h,x\in\Gamma.$$
 We thank N. Monod for suggesting us to look for irreducible affine actions of $\Gamma\times\Gamma$ with linear part $\vartheta$.

\begin{Prop}\label{amenICC}
Let $\Gamma$ be an infinite, countable, amenable ICC group. Then $\vartheta$ is the linear part of some irreducible affine action of $\Gamma\times\Gamma$, which can be chosen to have almost fixed points. 
\end{Prop}

\bpr Since $\Gamma$ is amenable and infinite, the representation $\vartheta$ almost has invariant vectors but no non-zero fixed vector.  Hence the space $B^1(\Gamma\times\Gamma,\vartheta)$ is not closed in  $Z^1(\Gamma\times\Gamma,\vartheta)$, by  \cite[Corollaire1]{Gui2} (note that countability is used here). Choose a cocycle $b$ in the closure of $B^1$ but not in $B^1$. Then the corresponding affine action $\alpha_{\vartheta,b}$ almost has fixed points. Finally, note that $\vartheta$ is an irreducible representation of $\Gamma\times\Gamma$, as $\Gamma$ is ICC. So $\alpha_{\vartheta,b}$ is irreducible, by Example \ref{basicex2}. \epr

\medskip
This must be contrasted with Theorem \ref{thom} above, which deals with the left regular representation of an amenable group.

\subsection{A super-rigidity result}

We now reach a super-rigidity result for lattices in a product of locally compact groups.

\begin{Thm}\label{superrig} Let $G=G_1\times\dots\times G_n$ be the product of non-trivial, compactly generated, locally compact groups, and let $\Gamma$ be a lattice in $G$, projecting densely to all factors. 
Assume that either $\Ga$ is co-compact, or that every $G_i$ is the group of $K_i$-points of an almost $K_i$-simple, $K_i$-isotropic linear algebraic group over some local field $K_i$.
Let $\pi$ be a unitary representation of $\Gamma$, not weakly containing the trivial representation, and let $\alpha$ be an affine isometric action of $\Gamma$ with linear part $\pi$. The following properties are equivalent:
\begin{enumerate}
\item[i)] $\alpha$ is irreducible.
\item[ii)]  For every $i=1,\dots,n$, there exists an irreducible affine action $\alpha_i$ of $G$, with $\alpha_i$ factoring through $G_i,$ such that $\alpha\simeq(\bigoplus_{i=1}^n \alpha_i)|_\Gamma$.
\end{enumerate}
\end{Thm} 

\bpr $(ii)\Rightarrow(i)$ follows by induction over $n$, combining Proposition \ref{irredproduct} with Theorem \ref{neretin} (and appealing to Remark~\ref{Rem-IntegrableCocycle} in the non-co-compact case).

$(i)\Rightarrow(ii)$ Let $b\in Z^1(\Gamma,\pi)$ be the cocycle defining $\alpha$. By a result of Shalom (\cite{ShaInv}, Corollary 4.2, using the assumption that $\pi$ does not weakly contain the trivial representation): $b$ is cohomologous to a sum $b_1+\dots+b_n$, where $b_i$ takes values in a $\pi(\Gamma)$-invariant subspace ${\cal H}_i\subset{\cal H}$; moreover, denoting by $\sigma_i$ the restriction of $\pi$ to ${\cal H}_i$, the affine action $\alpha_{\sigma_i,b_i}$ extends continuously to an affine action $\alpha_i$ of $G$ that factors through an action of $G_i$.

As in the proof of Proposition \ref{irredproduct}, conjugating $\alpha$ by a translation we may assume $b=b_1+\dots+b_n$, from which we deduce $\alpha=(\alpha_1\oplus\dots\oplus\alpha_n)|_\Gamma$. Since $\alpha_{\sigma_i,b_i}$ is a projected action of $\alpha$, it is an irreducible action of $\Gamma$. Finally, since $\alpha_i|_\Gamma=\alpha_{\sigma_i,b_i}$ and the projection of $\Gamma$ to $G_i$ is dense, $\alpha_i$ is an irreducible action of $G$.\epr

\begin{Cor}\label{irredH1} Keep notations as in Theorem \ref{superrig}. Let $\pi$ be a unitary irreducible representation of $\Gamma$, not containing weakly the trivial representation. If $H^1(\Gamma,\pi)\neq 0$, then for some $i=1,\dots,n$ the representation $\pi$ extends to a unitary irreducible representation $\sigma_i$ of $G$ factoring through $G_i$. Moreover the restriction map ${\rm Rest}_G^\Gamma: H^1(G,\sigma_i)\rightarrow H^1(\Gamma,\pi)$ is an isomorphism.
\end{Cor}

\bpr The first statement is obtained from Theorem \ref{superrig} exactly as the same way as Corollary \ref{irredproduct2} was obtained from Proposition \ref{irredproduct}. It also shows surjectivity of ${\rm Rest}_G^\Gamma$. Injectivity follows immediately from density of the projection of $\Gamma$ in $G_i$.\epr

\begin{Ex} \rm\begin{enumerate}
\item[i)] Let $p$ be a prime number. The group $PSL_2(\QQ_p)$ has a unique unitary irreducible representation $\sigma$ with non-vanishing $H^1$ (it is the representation on the first $L^2$-cohomology of the Bruhat-Tits tree); similarly $PSL_2(\RRR)$ has two unitary irreducible representations $\pi_+,\pi_-$ with non-vanishing $H^1$ (these are the representations on square-integrable holomorphic and anti-holomorphic 1-forms on the Poincar\'e disk); for all this, see \cite{BoWa}. Viewing $\Ga_p=:PSL_2(\mathbf{Z}[\frac{1}{p}])$ as a lattice in $PSL_2(\QQ_p)\times PSL_2(\RRR)$, we see from Corollary \ref{irredH1} that $\Ga_p$ has exactly three irreducible unitary representations, not weakly containing the trivial representation, with non-vanishing $H^1$, namely the restrictions of $\sigma,\pi_+,\pi_-$ to $\Ga_p$.

Similarly, viewing $\Lambda_p=:PSL_2(\mathbf{Z}[\sqrt{p}])$ as a lattice in $PSL_2(\RRR)\times PSL_2(\RRR)$, we see that $\Lambda_p$ has exactly four unitary irreducible representations, not weakly containing the trivial representation, with non-vanishing $H^1$: namely, $\pi_+|_{\Lambda_p},\,\pi_-|_{\Lambda_p},\,\pi_+\circ\tau,\,\pi_-\circ\tau$, where $\tau:a+b\sqrt{p}\mapsto a-b\sqrt{p}$ is the non-trivial element of the Galois group ${\rm Gal}(\QQ(\sqrt{p})/\QQ)$.
\item[ii)] Consider the quadratic form $Q$ in 5 variables, defined over $\QQ(\sqrt{2})$:
 $$Q(x)=x_1^2+x_2^2+x_3^2+\sqrt{2}x_4^2-x_5^2.$$
Set $\Ga=SO_0(Q)(\mathbf{Z}[\sqrt{2}])$, and view it as a lattice in $G=SO_0(Q)(\RRR)\times SO_0(\tau Q)(\RRR)$, where as above $\tau$ denotes the non-trivial element of the Galois group ${\rm Gal}(\QQ(\sqrt{2})/\QQ)$. As a Lie group $G$ is isomorphic to $SO_0(4,1)\times SO_0(3,2)$, the latter factor having property (T), the former not. Actually it is known (see \cite{BoWa}) that $SO_0(4,1)$ has a unique irreducible unitary representation $\pi$ with non-zero $H^1$. By Corollary \ref{irredH1}, the group $\Ga$ has a unique irreducible unitary representation, not weakly containing the trivial representation, with non-zero $H^1$: it is $\pi|_\Ga$.
\end{enumerate}
\end{Ex}

\section{On the first $L^2$-Betti number of a locally compact group}

Let $G$ be a unimodular, locally compact group with Haar measure $dg.$
Recall that a unitary irreducible representation $(\sigma, {\H}_\sigma)$ of $G$ is {\it square-integrable}
if
$$
\int_G |\langle \sigma(g) \xi|\xi\rangle |^2 dg <\infty \qquad \mbox{for all}\; \xi\in \H_\sigma.
$$
 This is the case if and only if    $\sigma$
 is a subrepresentation of the left regular representation $(\lambda_G, L^2(G))$
 of $G$.  Indeed, there exists a constant $d_\sigma>0,$ called the \emph{formal dimension}
 of $\sigma,$ such that the orthogonality relations hold
 $$
\int_G \langle \sigma(g) \xi|\eta\rangle\overline{\langle \sigma(g) \xi'|\eta'\rangle} dg = d_\sigma^{-1}\langle  \xi|\xi'\rangle\langle  
\overline{\eta|\eta'\rangle}
\qquad \mbox{for all}\, \xi,\xi', \eta, \eta'\in \H_\sigma
$$
 For every unit vector $\xi_0\in \H_\sigma,$ the $G$-equivariant  map
 $L:\mathcal{H}_\sigma \to  L^2(G)$ given by $L\eta (g) = \sqrt{d_\sigma}  \langle  \pi(g^{-1}) \eta|\xi_0\rangle,$
 is isometric and identifies $\H_\sigma$ with a $\la_G(G)$-invariant closed subspace of $L^2(G).$
  
 We denote by $\hat{G}_d$ the {\it discrete series} of $G$, i.e. the set of equivalence classes of square-integrable representations. Let $\Ga$ be a lattice in $G$.

Fix $\sigma\in \widehat G_d$ with  Hilbert space $\H_\sigma.$ The restriction of $\sigma$ to $\Ga$ extends to 
$L(\Ga)$ so that 
$\H_\sigma$ is a Hilbert module over  $L(\Ga).$  As such, $\H_\sigma$ has a 
 von Neumann dimension ${\rm dim}_{L(\Ga)} \H_\sigma$ (see Subsection~\ref{SS: LeftRegRep}).
 This dimension  is given by  Atiyah-Schmid's formula from \cite{AS} (see also \cite[Theorem 3.3.2]{GHJ}):
$$\dim_{L(\Gamma)}{\cal H}_\sigma =d_\sigma\textnormal{covol}(\Gamma).$$

As in  Subsection~\ref{SS: LeftRegRep}, set 
$$\beta^1_{(2)} (\Ga) ={\rm dim}_{R(\Ga)}  H^1_{(2)}(\Gamma),$$

\begin{Thm}\label{betti+discreteser+NonICC} 
Let  $G$ be  separable, compactly generated, locally compact group containing a finitely generated lattice $\Ga$  satisfying condition (\ref{funddomain}) from Remark~\ref{Rem-IntegrableCocycle} (e.g., $\Gamma$ co-compact).
Assume that $G$ is not amenable. 
Then 
$$\beta^1_{(2)}(\Ga)\geq \textnormal{covol}(\Gamma)\sum_{\sigma\in\hat{G}_d}  d_\sigma\cdot\dim_{\mathbf{C}} H^1(G,\sigma).$$
\end{Thm}

\bpr  
 It is enough to prove that, for every finite subset $F$ of $\hat{G}_d$ and  integers $k_\sigma$ 
 with $k_\sigma\leq\dim_{\mathbf{C}} H^1(G,\sigma)$ for $\sigma \in F,$ we have 
$$\beta^1_{(2)} (\Ga)\geq \textnormal{covol}(\Gamma)\sum_{\sigma\in F} k_\sigma   d_\sigma.$$ 
Choose 1-cocycles $b_1,\dots, b_{k_\sigma}$ such that the classes $[b_1],\dots, [b_{k_\sigma}]$ are linearly independent in $H^1(G,\sigma)$ and form the affine isometric action 
$$\alpha=\bigoplus_{\sigma\in F}(\oplus_{i=1}^{k_\sigma}\alpha_{\sigma,b_i});$$
Propositions \ref{linearlyindep} and \ref{directsum} implies that  the affine action $\alpha$ is irreducible.

By Theorem \ref{neretin}, the restriction $\alpha|_\Ga$ is irreducible. 
Moreover, $\Ga$ is non amenable as $G$ is non amenable.
 Hence, by Corollary~\ref{Cor1-AffRegRep} combined with the Atiyah-Schmid formula from above, we  have
\begin{equation}\label{discreteser}
\beta^1_{(2)}(\Ga)\geq \sum_{\sigma\in F} k_\sigma {\rm dim}_{L(\Ga)} {\H}_\sigma= 
 \textnormal{covol}(\Gamma)\sum_{\sigma\in F} k_\sigma   d_\sigma
\end{equation}
 \epr

\medskip
Let $G$ be a  second countable, locally compact unimodular group with Haar measure $dg.$
Denote by $L(G)$ the group von Neumann algebra of $G$; it carries a semi-finite trace $\psi$ defined on the positive cone of $L(G)$ by $\psi(x^*x)=\int_G |f(g)|^2\,dg$, where $x$ is left convolution by $f\in L^2(G)$; note that $\psi$ depends on the choice of the  Haar measure on $G$.

In two papers \cite{Pet, KPV}, Petersen and Kyed-Petersen-Vaes extended the classical definition of $L^2$-Betti numbers for discrete groups \cite{CG} to that more general framework, by setting
$$\beta^n_{(2)}(G):=\dim_{(L(G),\psi)} H^n(G,\lambda_G)$$
where $\lambda_G$ denotes the left regular representation on $L^2(G)$, and $\dim_{(L(G),\psi)}$ denotes the von Neumann dimension of $L(G)$-modules with respect to the semi-finite trace $\psi$.
They established a number of important results; in particular $\beta^1_{(2)}(G)<\infty$ as soon as $G$ is compactly generated, and 
$$\beta^n_{(2)}(G)=\frac{\beta^n_{(2)}(\Gamma)}{\textnormal{covol}(\Gamma)}$$ for every lattice $\Gamma$ in $G$.

 Recall that a locally compact group which contains a lattice is unimodular.
 
\begin{Thm}\label{betti+discreteser} Let $G$ be a second countable, compactly generated, locally compact group. Assume that 
 $G$ contains a finitely generated lattice satisfying condition (\ref{funddomain}) from Remark~\ref{Rem-IntegrableCocycle} (e.g. a co-compact lattice). Then
$$\beta^1_{(2)}(G)\geq \sum_{\sigma\in\hat{G}_d} d_\sigma\cdot\dim_{\mathbf{C}} H^1(G,\sigma).$$
\end{Thm}

\bpr  When  $G$ is not amenable,  the inequality is a direct consequence
of Theorem~\ref{betti+discreteser+NonICC}  and the formula linking $\beta^1_{(2)}(G)$ and $\beta^n_{(2)}(\Ga)$
from \cite{Pet, KPV}.

So we may assume that $G$ is amenable. We claim that both sides of the equality are zero. 
The vanishing of $\beta^1_{(2)}(G)$ follows from Theorem C in \cite{KPV}.

 Now we check that $H^1(G,\sigma)=0$ for every $\sigma\in\hat{G}_d$. 
  By (2.10) in \cite{KPV}, the vanishing of $\beta^1_{(2)}(G)$ implies 
  that the reduced first cohomology group $\overline{H}^1(G,\lambda_G)$
  is trivial.  
   Since $\sigma$ is a subrepresentation of 
  $\lambda_G$, we get $\overline{H}^1(G,\sigma)=0$. 
  
  Assume first that  $\sigma$ is not the trivial representation $1_G$. Since $\sigma$ is square-integrable, it defines a closed point in the dual $\hat{G}$.
  So, $\sigma$ does not weakly contain $1_G$ and hence $B^1(G,\sigma)$ is closed in $Z^1(G,\sigma),$
  by   \cite[Th\'eor\`eme 1]{Gui2};  therefore $H^1(G,\sigma)=0$.
  
  On the other hand, if $\sigma$ is the trivial representation $1_G$, then $G$ must be compact and therefore
  ${H}^1(G,1_G)={\rm Hom}(G,\mathbf{C})=0. $ \epr

\medskip
\noindent 
\begin{Rem}
\label{Rem-discreteser}
\rm
The proof  of  Theorem~\ref{betti+discreteser}   shows that the conclusion of 
Theorem~\ref{betti+discreteser+NonICC} holds also in the case where $G$ is amenable.
\end{Rem}

\begin{Cor} Keep the assumptions of Theorem \ref{betti+discreteser}.
 If $\beta^1_{(2)}(G)=0$, then $H^1(G,\sigma)= 0$ for all $\sigma\in\hat{G}_d$. \epr
\end{Cor}

\begin{Cor}\label{trees} Let $X_{k,\ell}$ be the $(k,\ell)$-biregular tree ($k=\ell$ being allowed). Let $G$ be a closed non-compact subgroup of ${\rm Aut}(X_{k,\ell})$, acting transitively on the boundary $\partial X_{k,\ell}$ and with two orbits on vertices of $X_{k,\ell}$. Normalize Haar measure on $G$ so that edge stabilizers have measure 1. Let $\sigma_0$ be the unique irreducible, square-integrable representation of $G$ with non-vanishing $H^1$ (see \cite{Neb}). Then $1-\frac{1}{k}-\frac{1}{\ell}\geq d_{\sigma_0}.$
\end{Cor}

\bpr First, $G$ contains a co-compact lattice (by Theorem 3.10 in \cite{BLub}), so we may apply Theorem \ref{betti+discreteser}:
$$\beta^1_{(2)}(G)\geq d_{\sigma_0}\dim_{\mathbf{C}}H^1(G,\sigma_0).$$
Second, $\beta1_{(2)}(G)=1-\frac{1}{k}-\frac{1}{\ell}$ by Corollary 5.18 in \cite{Pet}. Third, $\dim_{\mathbf{C}}H^1(G,\sigma_0)=1$ by the main Theorem in \cite{Neb}.
\epr

\medskip
\noindent 
\begin{Rem}
\label{Rem-discreteser2}
\rm
Theorem \ref{betti+discreteser} served as motivation for the main result in \cite{PetVal}: for $G$ a type I, unimodular, separable, locally compact group:
$$\beta^n_{(2)}(G)=\sum_{\sigma\in\hat{G}_d} d_\sigma\cdot\dim_{\mathbf{C}} H^n(G,\sigma)+\int_{\hat{G}\backslash\hat{G}_d}\dim_{\mathbf{C}}\overline{H}^n(G,\omega)\,d\mu(\omega),$$
where $\mu$ is the Plancherel measure on the dual $\hat{G}$ of $G$, and $\overline{H}^n$ denotes reduced $n$-cohomology. The proof is completely different. Observe the different sets of assumptions: type I in \cite{PetVal}, existence of a suitable lattice in Theorem \ref{betti+discreteser} above. 

For infinite discrete groups, Theorem \ref{betti+discreteser} just gives $\beta^1_{(2)}(G)\geq 0$, since $\hat{G}_d$ is empty in this case. On the other hand, the computations in \cite{PetVal} show that equality may occur either in Theorem \ref{betti+discreteser} or in Corollary \ref{trees}, with the right hand side being non-zero: this is the case for $PSL_2(\mathbf{R}),\,PSL_2(\mathbf{Q}_p)$ and for ${\rm Aut}(X_{k,\ell})$. Actually it follows from \cite{PetVal} that equality holds in Corollary \ref{trees}, under the extra assumption that $G$ is type I. It is an open question whether a group satisfying the assumptions of Corollary \ref{trees} must be type I.
\end{Rem}
\section{Comparison with other forms of irreducibility}

In \cite{CTV-SO}, the authors study orbits of affine isometric actions, and make the following definitions.

\begin{Def} \rm \begin{enumerate}[i)]
\item An affine isometric action $\alpha$ of a group $G$ on a Hilbert space ${\cal H}$ has enveloping orbits if the closed convex hull of every orbit is equal to ${\cal H}$.
\item A unitary representation $\pi$ of $G$, is strongly cohomological if $H^1(G,\sigma)\neq 0$ for every non-zero sub-representation $\sigma$ of $\pi$.
\end{enumerate}\end{Def}

It is observed in Lemma 4.3 of \cite{CTV-SO} that the linear part of an action with enveloping orbits is strongly cohomological. We notice that irreducibility lies in between having enveloping orbits and having a strongly cohomological linear part.

\begin{Prop}\label{enveloping} Let $\pi$ be a unitary representation of $G$. Every of following properties implies the next one:
\begin{enumerate}[i)]
\item There exists an affine isometric action with linear part $\pi$ and with enveloping orbits.
\item There exists an irreducible affine isometric action with linear part $\pi$.
\item $\pi$ is strongly cohomological.
\end{enumerate}
\end{Prop}

\bpr $(i)\Rightarrow(ii)$ It follows from the definitions that, if an affine isometric action has enveloping orbits, then it is irreducible.

$(ii)\Rightarrow(iii)$ This follows from $(A1)\Rightarrow(A3)$ in Proposition \ref{affine}.
\epr

Let us check that none of the converse implications  in Proposition \ref{enveloping} holds.

\begin{Ex} \rm Let $F=\mathbf{F}_2$ be the free group on 2 generators. As observed in Remark 3.5 of \cite{CTV-SO}, every representation of $G$ is strongly cohomological. Now let $\pi$ be the trivial representation of $G$ on a Hilbert space with dimension $>2$. There is no irreducible affine isometric action with linear part $\pi$.
\end{Ex}
The following example, suggested by Y. Cornulier, shows that the converse of $(i)\Rightarrow(ii)$ in Proposition \ref{enveloping} does not hold in infinite dimension. We denote by ${\rm Sym}(\mathbf{N})$ the full symmetric group of $\mathbf{N}$ (viewed as a discrete group), and by $C_2^{(\mathbf{N})}$ the direct sum of countably many copies of the cyclic group $C_2$ of order 2. Note that ${\rm Sym}(\mathbf{N})$ acts on $C_2^{(\mathbf{N})}$ by permutation of the indices.

\begin{Prop}\label{YC} Let $G$ be the semi-direct product $C_2^{(\mathbf{N})}\rtimes {\rm Sym}(\mathbf{N})$. Then $G$ admits an irreducible representation $(\pi,{\cal H})$ and an unbounded 1-cocycle $b\in Z^1(G,\pi)$ such that, for a dense set of vectors $w\in{\cal H}$, the function $g\mapsto\langle b(g)|w\rangle$ is bounded on $G$ (so that $\overline{{\rm conv}(b(G))}\neq{\cal H}$).
\end{Prop}

\bpr We identify $C_2$ with the multiplicative group $\{\pm 1\}$, and $C_2^{(\mathbf{N})}$ with the group of finitely supported functions $\mathbf{N}\rightarrow\{\pm 1\}$. Let ${\cal F}$ be the space of all real-valued sequences on $\mathbf{N}$, and ${\cal H}=\ell^2$ be the subspace of square-summable sequences. Then $C_2^{(\mathbf{N})}$ acts on ${\cal F}$ by pointwise multiplication, and ${\rm Sym}(\mathbf{N})$ acts on ${\cal F}$ by permutation of the indices. Let $\sigma$ be the corresponding linear representation of $G$ on ${\cal F}$. The subspace ${\cal H}$ is invariant, and we denote by $\pi$ the restriction of $\sigma$ to ${\cal H}$. The proof of the proposition will be carried out in four steps.
\begin{itemize}
\item[(i)] Clearly, the only $\sigma(G)$-fixed vector in ${\cal F}$ is $0$.
\item[(ii)] The representation $\pi$ is irreducible. Actually $\pi|_{{\rm Sym}(\mathbf{N})}$ is already irreducible. Indeed, by transitivity of the action of ${\rm Sym}(\mathbf{N})$ on $\mathbf{N}$, we can  identify (in a ${\rm Sym}(\mathbf{N})$-equivariant way) $\mathbf{N}$ with ${\rm Sym}(\mathbf{N})/{\rm Sym}(\mathbf{N})_0$, where ${\rm Sym}(\mathbf{N})_0$ is the stabilizer of $0$ in ${\rm Sym}(\mathbf{N})$. So $\pi$ is equivalent to the quasi-regular representation on $\ell^2({\rm Sym}(\mathbf{N})/{\rm Sym}(\mathbf{N})_0)$. Now observe that ${\rm Sym}(\mathbf{N})_0$ is equal to its commensurator in ${\rm Sym}(\mathbf{N})$; indeed, for $g\in {\rm Sym}(\mathbf{N})\backslash {\rm Sym}(\mathbf{N})_0$, the subgroup ${\rm Sym}(\mathbf{N})_0\cap g{\rm Sym}(\mathbf{N})_0g^{-1}$ is the stabilizer of $g(0)$ in ${\rm Sym}(\mathbf{N})_0$, so it has infinite index as ${\rm Sym}(\mathbf{N})_0$ acts transitively on $\mathbf{N}\backslash\{0\}$. Irreducibility then follows from Mackey's classical criterion 
for irreducibility of induced representations from self-commensurating subgroups \cite{Mac}.
\item[(iii)] Let $v=(1,1,1,\dots)$ be a constant sequence in ${\cal F}$. Form the affine action $t_{v}\circ\sigma\circ t_{-v}$. The associated 1-cocycle is $b(g)=v-\sigma(g)v$, which is $0$ if $g\in {\rm Sym}(\mathbf{N})$ and has finite support if $g\in C_2^{(\mathbf{N})}$. In particular, this affine action preserves ${\cal H}$ and induces on it an affine action $\alpha$. Since $v$ is the only fixed point of $t_{v}\circ\sigma\circ t_{-v}$ (as seen above) and $v\notin{\cal H}$, we see that $\alpha$ has no fixed point, i.e. $b$ is unbounded as a map $G\rightarrow{\cal H}$. Note also that $\alpha$ is irreducible, since $\pi$ is.
\item[(iv)] Observe that $b(G)$ is the set of sequences consisting of $0$'s and $2$'s, with finitely many $2$'s. View $\ell^1$ as a dense subspace of $\ell^2$. For $w\in\ell^1$ and $g\in G$, we have $|\langle b(g)|w\rangle|=|\sum_{n=0}^{\infty} b(g)_n w_n|\leq 2\|w\|_1$.
\end{itemize}\epr

\medskip
It turns out that, in Proposition \ref{enveloping}, the converse of $(ii)\Rightarrow(i)$ holds in {\it finite} dimension.

\begin{Prop} Let $\alpha$ be an affine isometric action of a group $G$ on $\mathbf{R}^n$. If $\alpha$ is irreducible, then $\alpha$ has enveloping orbits.
\end{Prop}

\bpr We first observe that the result trivially holds for $n=1$, since by irreducibility $\alpha(G)$ must contain a non-zero translation. Now, proceeding by contradiction, let $n$ be the smallest integer such that there exists an irreducible affine isometric action $\alpha$ on $\mathbf{R}^n$, with the property that for some orbit $\alpha(G)x_0$, the closed convex set $C:=\overline{\textnormal{conv}(\alpha(G)x_0)}$ is not equal to $\mathbf{R}^n$. Then $C$ is contained in some closed affine half-space $\{x\in\mathbf{R}^n:\langle x|w\rangle\leq a\}$, for some unit vector $w\in\mathbf{R}^n$ and some $a\in\mathbf{R}$. As $C$ is unbounded, it contains some half-line $D=x_0+\mathbf{R}^+.v_0$, where $v_0$ is some unit vector, such that $\langle w|v_0\rangle\leq 0$. Since $\alpha(g)D\subset C$ for every $g\in G$, we have similarly $\langle w|\pi(g)v_0\rangle\leq 0$ for every $g\in G$. Now two cases may occur:
\begin{itemize}
\item $\langle w|\pi(g)v_0\rangle< 0$ for some $g\in G$. Let $K$ be the closure of $\pi(G)$ in the orthogonal group $O(n)$. So $K$ is a compact group, with normalized Haar measure $dk$. Set $v=\int_K k.v_0\,dk$: then $v\neq 0$ since $\langle w|v\rangle=\int_K \langle w|k.v_0\rangle\,dk<0$ (as the integrand is $<0$ on a neighbourhood of $\pi(g)$). So $v$ is a non-zero $\pi(G)$-fixed vector. Let then $\alpha_0$ be the projected action on the 1-dimensional subspace $V=\mathbf{R}.v$; the action $\alpha_0$ is irreducible by Proposition \ref{affine}. On the other hand, the projection of $\alpha(G)x_0$ is contained in a half-line of $V$, contradicting the result for $n=1$. 
\item $\langle w|\pi(g)v_0\rangle=0$ for every $g\in G$. Let then $V_0$ be the $\pi(G)$-invariant subspace spanned by the $\pi(g)v_0$'s, let $V_1$ be the orthogonal of $V_0$, and let $\pi_1$ be the restriction of $\pi$ to $V_1$. Let $\alpha_1$ be the projected action on $V_1$. By Proposition \ref{affine}, $\alpha_1$ is irreducible, so it has enveloping orbits by minimality of $n$. On the other hand the projection of $\alpha(G)x_0$ on $V_1$ is contained in a closed affine half-space, a contradiction.
\end{itemize}
\epr

\vspace{5mm}

Bachir Bekka\\
 IRMAR UMR 5525\\
 Universit\'e de  Rennes 1\\
Campus Beaulieu\\
 F-35042  Rennes Cedex, France\\
\verb"bachir.bekka@univ-rennes1.fr"

\vspace{5mm}
Thibault Pillon and Alain Valette\\
Institut de Math\'ematiques\\
Universit\'e de  Neuch\^atel  \\
11 Rue Emile Argand\\
CH-2000 Neuch\^atel, Switzerland\\
\verb"thibault.pillon@unine.ch; alain.valette@unine.ch"

\end{document}